\documentclass[11pt]{article}

\usepackage{tikz}
\usetikzlibrary{calc}
\usetikzlibrary{matrix}
\usepackage{amsmath,amssymb,amsfonts,amsthm} 
\usepackage{verbatim}
\usepackage{hyperref}                 

\usepackage{graphicx}
\allowdisplaybreaks

\def \p{\partial}
\newcommand{\V}[1]{\mbox{\boldmath $ #1 $}}

\newcommand{\bey}{\begin{eqnarray}}
\newcommand{\eey}{\end{eqnarray}}
\newcommand{\nn}{\nonumber}
\newcommand{\beq}{\begin{equation}}
\newcommand{\eeq}{\end{equation}}
\theoremstyle{plain}
\newtheorem{thm}{\hspace{6mm}Theorem}[section]

\newtheorem{lem}{\hspace{6mm}Lemma}[section]

\theoremstyle{definition}

\theoremstyle{remark}
\newtheorem{exam}{\hspace{6mm}Example}[section]
\newtheorem{rem}{\hspace{6mm}Remark}[section]

\vspace{5mm}

\begin{document}

\date{}
\title{Maximum principle in linear finite element approximations of anisotropic diffusion-convection-reaction problems
}
\author{Changna Lu\thanks{College of Math and Statistics,
Nanjing University of Information Science and Technology,
Nanjing, Jiangsu 210044, China. ({\tt luchangna@nuist.edu.cn})},
\and Weizhang Huang\thanks{
Department of Mathematics, the University of Kansas, Lawrence, KS 66045,
U.S.A. ({\tt huang@math.ku.edu})},
\and
Jianxian~Qiu\thanks{School of Mathematical Sciences, Xiamen University, Xiamen,
Fujian 361005, China. ({\tt jxqiu@xmu.edu.cn})}
}
\maketitle

\vspace{10pt}

\begin{abstract}
A mesh condition is developed for linear finite element approximations of anisotropic
diffusion-convection-reaction problems to satisfy a discrete maximum principle.
Loosely speaking, the condition requires that the mesh be simplicial and $\mathcal{O}(\|\V{b}\|_\infty h
+ \|c\|_\infty h^2)$-nonobtuse when the dihedral angles are measured in the metric specified by the inverse
of the diffusion matrix, where $h$ denotes the mesh size and $\V{b}$ and $c$ are the coefficients
of the convection and reaction terms. In two dimensions, the condition can be replaced by a weaker mesh condition
(an $\mathcal{O}(\|\V{b}\|_\infty h + \|c\|_\infty h^2)$ perturbation of a generalized Delaunay condition).
These results include many existing mesh conditions as special cases. Numerical results are presented
to verify the theoretical findings.
\end{abstract}



\vspace{10pt}

\section{Introduction}

We are concerned with the linear finite element (FEM) solution of the anisotropic diffusion equation
\beq
\label{bvp-pde}
 -\nabla \cdot (\mathbb{D} \, \nabla u) + \V { b } \cdot \nabla u + c \, u = f,\qquad \mbox{ in } \quad \Omega
\eeq
subject to the Dirichlet boundary condition
\beq
\label{bvp-bc}
u =  g, \qquad \mbox{ on } \quad \partial \Omega
\eeq
where $\Omega \subset \mathbb{R}^d$ ($d \ge 1$) is a connected polyhedron and $\mathbb{D} = \mathbb{D}( \V{x} )
\in \mathbb{R}^{d \times d}$ (the diffusion matrix), $\V{b} = \V{b}(\V{x}) \in \mathbb{R}^d$, $c=c(\V{x})$, $f= f(\V{x})$,
and $g= g(\V{x})$ are given, sufficiently smooth functions defined on $\Omega$.
We assume that for any $\V{x}\in \Omega$, $\mathbb{D}( \V{x} )$ is symmetric and strictly positive definite
and functions $\V{b}$ and $c$ satisfy
\beq
 c(\V{x}) - \frac{1}{2} \nabla \cdot \V{b}(\V{x}) \geq 0 ,\qquad c (\V{x}) \geq 0 ,\qquad \forall \V{x}\in \Omega .
 \label{bvp-sc}
\eeq
It is known (e.g., see \cite{Eva98}) that the solution of the boundary value problem (BVP) (\ref{bvp-pde}) and
(\ref{bvp-bc}) satisfies the maximum principle.

The numerical solution of BVP (\ref{bvp-pde}) and (\ref{bvp-bc}) has attracted considerable attention from scientists and
engineers. The BVP is a prototype model for anisotropic diffusion problems which arise in various fields such as 
plasma physics \cite{GL09,GYK05,SH07}, petroleum reservoir simulation \cite{ABBM98a,MD06},
and image processing \cite{PM90,Wei98}.
Moreover, it has been amply demonstrated that a standard numerical method, such as a finite element, a finite difference,
or a finite volume method, does not necessarily satisfy a discrete maximum principle (DMP) and may produce
unphysical solutions that typically contain spurious oscillations, undershoots, and overshoots.
Furthermore, designing a numerical scheme to preserve the maximum principle is an important research
topic in its own right.
As a matter of fact, considerable work has been done in the past to develop numerical schemes to satisfy DMP;
e.g., see \cite{BKK08,BE04,Cia70,CR73,KK09,KKK07,KL95,Let92,Sto86,SF73,WaZh11,XZ99}
for isotropic diffusion problems ($\mathbb{D} = \alpha(\V{x}) I$ with $\alpha(\V{x})$ being a scalar function) and
\cite{DDS04,GL09,GYK05,Hua10,KSS09,LePot09,LH10,LSS07,LSSV07,LS08,MD06,SH07}
for anisotropic diffusion problems. In particular, it is shown in \cite{CR73} that the linear FEM satisfies DMP
when the mesh is simplicial and satisfies the so-called non-obtuse angle condition which requires that
the dihedral angles of all mesh elements be non-obtuse. In two dimensions the condition can be replaced
by a weaker condition (the Delaunay condition) which requires that the sum of any pair of angles opposite
a common edge is less than or equal to $\pi$ \cite{SF73}. Similar results have been obtained recently
for anisotropic diffusion problems in \cite{Hua10,LH10}.

It is pointed out that most of the existing work has been concerned with problems without convection terms.
For continuous problems, it is known (e.g., see \cite{Eva98}) that convection terms have no effect on
the satisfaction of the maximum principle by the solution. However, the situation is different for discrete
schemes. The main difficulty comes from the fact that discrete convection terms typically do not vanish
at an interior maximum point and the entries of the corresponding matrix can be both positive and negative.
A few researchers have tried to address the issue for isotropic diffusion problems.
For example, Xu and Zikatanov \cite{XZ99} employ a special number treatment for convection terms
so that they have no effect on the DMP satisfaction by the discrete solution.
Burman and Ern \cite{BE04} propose a nonlinear stabilized Galerkin approximation
of the Laplace operator which satisfies DMP on arbitrary meshes and for arbitrary space
dimension without resorting to the non-obtuse angle condition. They prove that the result can extend to
diffusion-convection-reaction problems with constant diffusion coefficient when the mesh is locally
quasi-uniform. More recently, Wang and Zhang \cite{WaZh11} study quasilinear
isotropic diffusion-convection-reaction problems and show that linear finite element approximations
satisfy DMP when the mesh is $\mathcal{O}(\|\V{b}\|_\infty h + \|c\|_\infty h^2)$-acute
(i.e., the dihedral angles of all mesh elements
are less than or equal to $\frac{\pi}{2}-\gamma_1 \|\V{b}\|_\infty h -  \gamma_2\|c\|_\infty h^2 $
for some positive constants $\gamma_1$ and $\gamma_2$).
On the other hand, no work has been done for anisotropic diffusion-convection-reaction problems.

The objective of this paper is to develop a mesh condition for linear finite element approximations
of anisotropic diffusion-convection-reaction problems (\ref{bvp-pde}) and (\ref{bvp-bc}) in any dimension
to satisfy a discrete maximum principle. We shall use the approach of \cite{LH10} to show
the stiffness matrix associated with the linear finite element discretization to be an $M$-matrix
and have non-negative row sums, with the focus on the treatments of the convection and reaction terms.
We shall also investigate the two dimensional case where a weaker sufficient condition can be
developed.

The paper is organized as follows. 
A linear finite element discretization for BVP (\ref{bvp-pde}) and (\ref{bvp-bc}) is introduced in Section 2.
In Section 3 geometric properties of the gradient of linear basis functions are studied.
A general mesh condition valid in any dimension and a specific and weaker condition in two dimensions
are developed in Section 4, followed by numerical results in Section 5.
Finally, Section 6 contains conclusions.

\section{Linear finite element formulation}

We consider the linear finite element solution of BVP (\ref{bvp-pde}) and (\ref{bvp-bc}).
Assume that an affine family of simplicial meshes $\{ \mathcal{T}_h \}$ is given for $\Omega$. Let
\[
U_g = \{ v \in H^1(\Omega) \; | \: v|_{\p \Omega} = g\}.
\]
Denote by $U_{g^h}^h$ the linear finite element space associated with mesh $\mathcal{T}_h$,
where $g^h$ is a piecewise linear approximation to $g$ on the boundary.
A linear finite element approximation $u^h \in U_{g^h}^h$ to BVP (\ref{bvp-pde})
and (\ref{bvp-bc}) is defined by
\bey
&& \int_{\Omega} (\nabla v^h)^{T} \, \mathbb{D} \,
\nabla u^h d \V{x} + \int_{\Omega} v^h \; (\V{b} \; \cdot \; \nabla u^h) d \V{x} + \int_{\Omega} c \; u^h \, v^h d \V{x}
\nn \\
&& \qquad \qquad = \int_{\Omega} f \, v^h d \V{x} , \quad \forall \; v^h \in U_0^h .
\label{disc-0}
\eey
The above equation can be rewritten as
\bey
&& \sum_{K \in \mathcal{T}_h} |K|\; (\nabla v^h)^{T} \, \mathbb{D}_K \,
\nabla u^h + \sum_{K \in \mathcal{T}_h} \int_{K} v^h \; (\V{b} \cdot \nabla u^h) d \V{x}
\nn \\
&& \qquad \qquad
+ \sum_{K \in \mathcal{T}_h} \int_{K} c \; u^h \, v^h d \V{x}
= \sum_{K \in \mathcal{T}_h} \int_{K} f \, v^h d \V{x} ,
\quad \forall \; v^h \in U_0^h
\label{disc-1}
\eey
where $|K|$ is the volume of element $K$ and $\mathbb{D}_K$ is the integral average of $\mathbb{D}$ over $K$, viz.,
\beq
\label{def-D}
\mathbb{D}_K = \frac{1}{|K|} \int_K \mathbb{D} \; d\V{x} .
\eeq

Scheme (\ref{disc-1}) can be expressed in a matrix form.
Denote the numbers of the elements, vertices, and interior vertices of mesh $\mathcal{T}_h$
by $N$, $N_v$, and $N_{vi}$, respectively. Assume that the vertices are ordered in such a way that
the first $N_{vi}$ vertices are the interior vertices. Then $U_0^h$ and $u^h$ can be expressed as
\beq
U_0^h = \text{span} \{ \phi_1, \cdots, \phi_{N_{vi}} \} ,
\eeq
\beq
\label{soln-approx}
u^h = \sum_{j=1}^{N_{vi}} u_j \phi_j + \sum_{j=N_{vi}+1}^{N_{v}} u_j \phi_j ,
\eeq
where $\phi_j$ denotes the linear basis function associated with the $j$-th vertex, $\V{a}_j$.
The boundary condition (\ref{bvp-bc}) is approximated by
\beq
u_j = g(\V{a}_j), \quad j = N_{vi}+1, ..., N_v .
\label{fem-bc}
\eeq
Substituting (\ref{soln-approx}) into (\ref{disc-1}), taking $v^h = \phi_j$ ($j=1, ..., N_{vi}$), and
combining the resulting equations with (\ref{fem-bc}), we obtain the linear algebraic system
\beq
\label{fem-linsys}
A \, \V{u} = \V{f},
\eeq
where $\V{u} = (u_1,..., u_{N_{vi}}, u_{N_{vi}+1},..., u_{N_v})^T$,
$\V{f} = ( f_1, ..., f_{N_{vi}}, g_{N_{vi}+1}, ..., g_{N_v} )^T$,
\beq
A = \left [\begin{array}{cc} A_{11} & A_{12} \\ 0 & I \end{array} \right ],
\label{fem-matrix}
\eeq
and $I$ is the identity matrix of size $(N_{v} - N_{vi})$.
The entries of the stiffness matrix $A$ and the right-hand-side
vector $\V{f}$ are given by
\begin{align}
a_{ij} & = \sum_{K \in \mathcal{T}_h} |K| (\nabla \phi_i)^{T} \, \mathbb{D}_K \,
\nabla \phi_j + \sum_{K \in \mathcal{T}_h} \int_{K} \phi_i \; (\V{b} \cdot \nabla \phi_j) d \V{x}
\notag \\
& \qquad \qquad + \sum_{K \in \mathcal{T}_h} \int_{K} c \; \phi_j \, \phi_i d \V{x} ,
\quad i=1, ..., N_{vi},\; j=1, ..., N_{v}
\label{stiffness-matrix}
\\
f_i & = \sum_{K \in \mathcal{T}_h} \int_{K} f \, \phi_i d \V{x},\quad \mbox{ \hfill} i=1, ..., N_{vi}.
\label{rhs-1}
\end{align}
In the following sections we shall investigate under what condition on the mesh the solution of (\ref{disc-1}) satisfies
a maximum principle. A key to this investigation is to understand geometric properties of the gradient
of linear basis functions which are to be described in the next section.

\section{Geometric properties of the gradient of linear basis functions}

Let $K$ be an arbitrary simplex with vertices
$\V{a}_1, \V{a}_2,..., \V{a}_{d+1}$. Denote the face opposite to vertex $\V{a}_i$
(i.e. the face not having $\V{a}_i$ as its vertex) by $S_i$ and its unit inward (pointing to $\V{a}_i$)
normal by $\V{n}_i$. The distance (or height) from vertex $\V{a}_i$ to face $S_i$
is denoted by $h_i$. The result of the following lemma exists in literature; e.g., see \cite{BKK08,KL95,XZ99}.

\begin{lem}
\label{lem3.1}
For any simplex $K \in \mathbb{R}^d$, the gradient of linear basis function $\phi_i$ associated any vertex
$\V{a}_i$ ($i = 1, ..., d+1$) is given by
\beq
\nabla \phi_i = \frac{\V{n}_i}{h_i}.
\label{lem-1}
\eeq
\end{lem}

It is remarked that Brandts et al. \cite{BKK08} have obtained the same result using the so-called $\V{q}$-vectors
defined through the edge matrix of elements. Specifically, they show that
$\V{q}_i$, a $\V{q}$-vector associated with face $S_i$, is an inward normal to $S_i$, has the length
$1/h_i$, and is equal to $\nabla \phi_i$; i.e.,
\beq
\V{q}_i = \frac{1}{h_i} \V{n}_i = \nabla \phi_i, \quad i = 1, ..., d+1.
\label{bkk08-1}
\eeq
These $\V{q}$-vectors will be used frequently in the remaining of the paper.

\begin{figure}[t]
\centering
\begin{tikzpicture}[scale = 1]
\draw [thick] (-1,0) -- (3,0) -- (2, 2) -- cycle;
\draw [below] (-1,0) node {$\V{a}_3$};
\draw [below] (3,0) node {$\V{a}_1$};
\draw [right] (2,2) node {$\V{a}_2$};
\draw [->] (-0.6,0) arc (0:30:0.42);
\draw [right] (-0.5,0.25) node {$\alpha_{12}$};
\draw [<-] (2.75,0.5) arc (135:180:0.7);
\draw [left] (2.5,0.3) node {$\alpha_{23}$};
\draw [->] (1,-0.4) -- (1,0.4);
\draw [below] (1,-0.4) node {$\V{q}_2$};
\draw [->] (3.0,1.25) -- (2,0.75);
\draw [right] (3.0,1.25) node {$\V{q}_3$};
\draw [->] (0.1,1.6) -- (0.82,0.52);
\draw [right] (0.1,1.6) node {$\V{q}_1$};
\draw [-] (3,0) -- (1.77,1.846);
\draw [left] (2.1,1.2) node {$h_1$};
\end{tikzpicture}
\caption{A sketch of unit inward normals, dihedral angles, and heights of element $K$.}
\label{f1}
\end{figure}
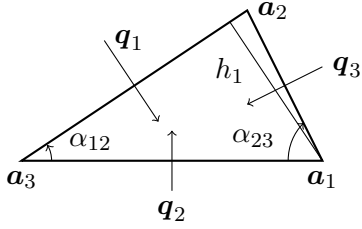

The next property of gradient of linear basis functions is related to the diffusion term
in stiffness matrix (\ref{stiffness-matrix}) for the case $\mathbb{D}_K = I$.
Denote the dihedral angle between any two faces $S_i$ and $S_j$
($i \neq j$) by $\alpha_{ij}$. It can be calculated as the supplement of the angle between
the inward normals to the faces, i.e.,
\beq
\label{dihedral}
\cos(\alpha_{ij}) = - \;\V{n}_i \cdot \V{n}_j = - \; \frac{\V{q}_i \cdot \V{q}_j}{\|\V{q}_i\| \cdot \|\V{q}_j\|} ,
\quad i \ne j.
\eeq
(In fact, (\ref{dihedral}) is often used as the definition of the dihedral angle.)
A sketch of the $\V{q}$-vectors, dihedral angles, and heights of an element are shown in Fig.~\ref{f1}.

The result of the following lemma is also known in literature; e.g., see \cite{BKK08,For89,Hua10}.

\begin{lem}
\label{lem3.2}
For any simplex $K \in \mathbb{R}^d$, we have
\beq
|K| (\nabla \phi_i )^T \nabla \phi_j = -\frac{|K|}{h_i h_j} \cos (\alpha_{ij} ), \quad i \neq j .
\label{lem3.2-1}
\eeq
It reduces to
\beq
|K| (\nabla \phi_i )^T \nabla \phi_j = -\frac{1}{2} \cot (\alpha_{ij} ), \quad i \neq j
\label{lem3.2-2}
\eeq
in two dimensions.
\end{lem}

\begin{proof}
Equation (\ref{lem3.2-1}) follows from Lemma \ref{lem3.1} and (\ref{dihedral}).

In two dimensions, $K$ is a triangle. Consider the case with $i=1$ and $j = 2$.
From Fig. \ref{f1}, we have
\[
|K| = \frac{h_1}{2} \|\V{a}_2 - \V{a}_3 \| = \frac{h_1 h_2}{2 \sin(\alpha_{12})} .
\]
Combining this result and (\ref{lem3.2-1}) gives (\ref{lem3.2-2}).
\end{proof}


We now study the diffusion term $|K| (\nabla \phi_i)^T \mathbb{D}_K \nabla \phi_j$
for general symmetric and positive definite matrix $\mathbb{D}_K$ using
Lemma \ref{lem3.2}.
Define
\beq
G_K(\V{x}) = \mathbb{D}_K^{-\frac 1 2} \V{x}:\;\; K \to \widetilde{K},
\label{G-1}
\eeq
where $\widetilde{K} = G(K)$. Obviously, $\widetilde{K}$ is also a simplex in $\mathbb{R}^d$.
For any vertex $\V{a}_i$, we denote the corresponding vertex, face, height, and $\V{q}$-vector of $\widetilde{K}$
by $\widetilde{\V{a}}_i$, $\widetilde{S}_i$, $\widetilde{h}_i$, and $\widetilde{\V{q}}_i$, respectively.
We have
\beq
\begin{cases}
& \widetilde{\V{a}}_i = \mathbb{D}_K^{-\frac 1 2} \V{a}_i, \quad
\widetilde{S}_i = \mathbb{D}_K^{-\frac 1 2} S_i,
|\widetilde{K}| = \mbox{det}(\mathbb{D}_K)^{-\frac 1 2} |K|,\\
& \widetilde{\V{q}}_i = \mathbb{D}_K^{\frac 1 2} \V{q}_i , \quad
\widetilde{h}_i = \| \V{q}_i\|_{\mathbb{D}_K}^{-1} ,
\end{cases}
\label{relation-1}
\eeq
where $\| \cdot \|_{\mathbb{D}_K}$ denotes the distance measured in the metric $\mathbb{D}_K$.
The derivations of the first three relations are trivial. To derive the last two, we first notice that
\[
\phi_i(\V{x}) = \phi_i(\mathbb{D}_K^{\frac 1 2} \widetilde{\V{x}}) = \widetilde{\phi}_i (\widetilde{\V{x}}) .
\]
Then from (\ref{bkk08-1}) we have
\[
\widetilde{\V{q}}_i = \widetilde{\nabla} \widetilde{\phi}_i = \mathbb{D}_K^{\frac 1 2} \nabla \phi_i
= \mathbb{D}_K^{\frac 1 2} \V{q}_i,
\]
which gives the second last relation in (\ref{relation-1}). The last relation is obtained by taking the norm of the
above equation.

To obtain the relation between $h_i$ and $\widetilde{h}_i$, we rewrite the last relation in (\ref{relation-1}) as
\[
\widetilde{h}_i = \frac{1}{\sqrt{(\V{q}_i)^T \mathbb{D}_K \V{q}_i}},
\]
from which we obtain
\beq
\frac{h_i}{\sqrt{\lambda_{max}(\mathbb{D}_K)}}\le \widetilde{h}_i \le \frac{h_i}{\sqrt{\lambda_{min}(\mathbb{D}_K)}},
\label{relation-2}
\eeq
where $\lambda_{max}(\mathbb{D}_K)$ and $\lambda_{min}(\mathbb{D}_K)$ denote the maximum and minimum
eigenvalues of $\mathbb{D}_K$, respectively.

Denote the dihedral angle between faces $\widetilde{S}_i$ and $\widetilde{S}_j$ by $\alpha_{ij,\mathbb{D}_K^{-1}}$.
Since $\widetilde{S}_i = \mathbb{D}_K^{-\frac 1 2} S_i $ and $\widetilde{S}_j=\mathbb{D}_K^{-\frac 1 2} S_j$,
it can also be viewed as the dihedral angle between $S_i$ and $S_j$
measured in the metric $\mathbb{D}_K^{-1}$.
Moreover, from \eqref{dihedral} we see that the angle can be calculated by
\beq
\cos (\alpha_{ij, \mathbb{D}_K^{-1}}) = - \frac{\widetilde{\V{q}}_i \cdot \widetilde{\V{q}}_j}
{\| \widetilde{\V{q}}_i \| \cdot \|\widetilde{\V{q}}_j \|}
= - \frac{\V{q}_i^T \mathbb{D}_K \V{q}_j}{\| \V{q}_i\|_{\mathbb{D}_K} \| \V{q}_j\|_{\mathbb{D}_K} } .
\label{dihedral-2}
\eeq

We now go back to the quantity $|K| (\nabla \phi_i)^T \mathbb{D}_K \nabla \phi_j$. Notice that
\[
|K| (\nabla \phi_i)^T \mathbb{D}_K \nabla \phi_j = \mbox{det}(\mathbb{D}_K)^{\frac 1 2}
|\widetilde{K}| (\widetilde{\V{q}}_i)^T \widetilde{\V{q}}_j.
\]
Applying Lemma \ref{lem3.2} to $\widetilde{K}$ and using relations (\ref{relation-1}), we have the following lemma.

\begin{lem}
\label{lem3.3}
For any simplex $K \in \mathbb{R}^d$ and any symmetric and positive definite matrix $\mathbb{D}_K$, we have
\beq
|K| (\nabla \phi_i )^T \mathbb{D}_K \nabla \phi_j = -\frac{|\widetilde{K}| \; \mbox{\em det}(\mathbb{D}_K)^{\frac 1 2}}
{\widetilde{h}_i \widetilde{h}_j} \cos (\alpha_{ij,\mathbb{D}_K^{-1}} ), \quad i \neq j .
\label{lem3.3-1}
\eeq
It reduces to
\beq
|K| (\nabla \phi_i )^T \mathbb{D}_K \nabla \phi_j = -\frac{\mbox{\em det}(\mathbb{D}_K)^{\frac 1 2}}{2}
\cot (\alpha_{ij,\mathbb{D}_K^{-1}} ), \quad i \neq j
\label{lem3.3-2}
\eeq
in two dimensions.
\end{lem}

\section{Mesh conditions for DMP satisfaction}

In this section we study the mesh conditions under which the linear finite element scheme (\ref{disc-1})
satisfies DMP. The main conclusions are given in Theorems \ref{thm4.1} and \ref{thm4.2}.

\begin{thm}
\label{thm4.1}
If the mesh satisfies
\begin{align}
& \frac{h_i^K}{\lambda_{min}(\mathbb{D}_K)} \cdot \frac{\|\V{b}\|_{\infty, K} }{(d+1)}
+ \frac{h_i^K h_j^K}{\lambda_{min}(\mathbb{D}_K)}\cdot \frac{\|c\|_{\infty, K} }{(d+1)(d+2)}
\le \cos (\alpha_{ij,\mathbb{D}_K^{-1}} ),
\label{thm4.1-1} \\
& \qquad \qquad \qquad \qquad i,j=1,...,d+1, \; i\neq j, \; \forall \; K \in \mathcal{T}_h 
\notag
\end{align}
where $\|\V{b}\|_{\infty, K} = \max_{\V{x}\in K} \| \V{b}(\V{x})\|$, $\|c\|_{\infty, K} = \max_{\V{x}\in K} c(\V{x})$,
and $h_i^K$'s and $\alpha_{ij,\mathbb{D}_K^{-1}}$'s are the heights and dihedral angles of element $K$,
respectively, then the linear finite element scheme (\ref{disc-1}) for BVP (\ref{bvp-pde}) and (\ref{bvp-bc})
satisfies DMP.
\end{thm}

\begin{proof}
Following \cite{LH10} we prove this theorem by showing that stiffness matrix $A$
defined in (\ref{fem-matrix}) and (\ref{stiffness-matrix}) has non-negative row sums and is
an M-matrix\footnote{Matrix $A$ is called an M-matrix if it is a Z-matrix (see (\ref{thm4.1-3})
and (\ref{thm4.1-4}) below) and satisfies $A^{-1} \ge 0$ (i.e., all entries of its inverse are nonnegative).}.
From Stoyan \cite[Theorem 1]{Sto86}, this implies that scheme (\ref{disc-1}) satisfies DMP.

(1) We first show that matrix $A$ has non-negative row sums. Notice that we only need to show that the first
$N_{vi}$ row sums are non-negative. Using the fact $\sum_{j=1}^{N_v}\phi_j( \V{x} )= 1$ and the assumption $c \geq 0$
(cf. (\ref{bvp-sc})), from (\ref{stiffness-matrix}) we have, for $i=1,...,N_{vi}$,
\bey
\sum_{j=1}^{N_v} a_{ij} &=& \sum_{K \in \mathcal{T}_h} |K| \;(\nabla \phi_i)^T \; \mathbb{D}_K \; \nabla \left( \sum_{j=1}^{N_v} \phi_j \right)
\nn \\
&&
+ \sum_{K \in \mathcal{T}_h} \;\int_K
\; \phi_i \; \left(\V{b} \cdot \nabla \left ( \sum_{j=1}^{N_v} \phi_j \right )\right) d \V{x}
\nn \\
&&
+ \sum_{K \in \mathcal{T}_h} \; \int_K
\; c \; \phi_i \; \left( \sum_{j=1}^{N_v} \phi_j \right) d \V{x}
\nn \\
&=& \sum_{K \in \mathcal{T}_h} \; \int_K
\; c \; \phi_i d \V{x}
\nn \\
&\geq& 0 .
\eey

(2) Next, we show that $A$ is a Z-matrix; i.e.,
\bey
a_{ij} &\leq& 0, \;\; \forall\; i \neq j, \; i = 1, ..., N_{vi},\; j = 1,..., N_{v}
\label{thm4.1-3}
\\
a_{ii} &\geq& 0, \;\; i = 1,..., N_{vi} .
\label{thm4.1-4}
\eey
Recall from Ciarlet \cite[Page 201]{Cia78} that
\beq
\int\limits_{K \in \omega_i} \phi_i d \V{x} = \frac{|K|}{d+1},\quad
\int\limits_{K \in \omega_i \cap \omega_j} \phi_i \phi_j d \V{x} = \frac{|K|}{(d+1)(d+2)} ,
\label{ciarlet78-1}
\eeq
where $\omega_i$ and $\omega_j$ are the element patches associated with vertices $\V{a}_i$ and
$\V{a}_j$, respectively. We havex
\begin{align}
a_{ij} &= \sum_{K \in \omega_i \cap \omega_j} \left (
|K| \;(\nabla \phi_i)^T \; \mathbb{D}_K \; \nabla \phi_j
+ \int_K \; \phi_i \; (\V{b} \cdot \nabla \phi_j ) d \V{x}
+ \; \int_K
\; c \; \phi_i \; \phi_j d \V{x} \right )
\tag{from (\ref{stiffness-matrix})} \\
&\leq \sum_{K \in \omega_i \cap \omega_j} \left ( |K| \; \;(\nabla \phi_i)^T \; \mathbb{D}_K \; \nabla \phi_j
+ \frac{1}{h_j^K} \int_K \phi_i | \V{b} \cdot \V{n}_j^K | d \V{x}
+ \int_K c \; \phi_i \; \phi_j \; d \V{x} \right )
\tag{Lemma \ref{lem3.1}} \\
&\leq \sum_{K \in \omega_i \cap \omega_j} \left ( |K| \; \;(\nabla \phi_i)^T \; \mathbb{D}_K \; \nabla \phi_j
+ \frac{|K| \; \|\V{b}\|_{\infty, K} }{h_j^K(d+1)}
+ \frac{|K|\; \|c\|_{\infty, K}}{(d+1)(d+2)} \right )
\tag{from (\ref{ciarlet78-1})} \\
&= \sum_{K \in \omega_i \cap \omega_j} \left (
-\frac{|K|}{\widetilde{h}_i^K \widetilde{h}_j^K} \cos (\alpha_{ij,\mathbb{D}_K^{-1}} )
+ \frac{|K| \; \|\V{b}\|_{\infty, K} }{h_j^K (d+1)}
+ \frac{|K| \; \|c\|_{\infty, K}}{(d+1)(d+2)} \right )
\tag{Lemma \ref{lem3.3}} \\
&= \sum_{K \in \omega_i \cap \omega_j} \frac{|K|}{\widetilde{h}_i^K \widetilde{h}_j^K} \left (
- \cos (\alpha_{ij,\mathbb{D}_K^{-1}} )
+ \frac{\widetilde{h}_i^K \widetilde{h}_j^K \|\V{b}\|_{\infty, K} }{h_j^K(d+1)}
+ \frac{\widetilde{h}_i^K \widetilde{h}_j^K\|c\|_{\infty, K}}{(d+1)(d+2)} \right )
\nn \\
&\le \sum_{K \in \omega_i \cap \omega_j} \frac{|K|}{\widetilde{h}_i^K \widetilde{h}_j^K} \left (
- \cos (\alpha_{ij,\mathbb{D}_K^{-1}} )
+ \frac{h_i^K \;\|\V{b}\|_{\infty, K} }{\lambda_{min}(\mathbb{D}_K) (d+1)}
\frac{}{} \right.
\nn \\
& \qquad \qquad \left. \frac{}{}
+ \frac{h_i^K h_j^K \; \|c\|_{\infty, K} }{\lambda_{min}(\mathbb{D}_K) (d+1)(d+2) } \right ) .
\tag{from (\ref{relation-2})}
\end{align}
Combining this with (\ref{thm4.1-1}) implies (\ref{thm4.1-3}).

On the other hand, for $i = 1,...,N_{vi}$,
\begin{align}
a_{ii} &= \sum_{K \in \mathcal{T}_h} |K| \;(\nabla \phi_i)^T \; \mathbb{D}_K \; \nabla \phi_i
+ \int_\Omega \; \phi_i \; (\V{b} \cdot \nabla \phi_i ) d \V{x}
+ \int_\Omega\; c \; \phi_i^2 d \V{x}
\tag{from (\ref{stiffness-matrix})} \\
& \ge \int_{\Omega} \phi_i ( \V{b} \cdot \nabla \phi_i ) d \V{x} + \int_{\Omega} c \; \phi_i^2 d \V{x}
\nn \\
& = \int_{\Omega} ( c - \frac{1}{2} \nabla \cdot \V{b} ) \phi_i^2 d \V{x} .
\tag{Gauss' divergence thm}
\end{align}
The assumption (\ref{bvp-sc}) implies that $a_{ii}\ge 0$. Thus, the stiffness matrix $A$ is a Z-matrix.

(3) We now show that $A_{11}$, the northwest block of matrix $A$, is an M-matrix. This is done by
showing $A_{11}$ is positive definite. For any vector $\V{v}=(v_1, v_2, ..., v_{N_{vi}})^T$, we define $v^h = \sum _{i=1}^{N_{vi}}v_i \phi_i \in U_0$. Notice that $\nabla v^h$ is constant on $K$.
As in the proof for $a_{ii}\ge 0$, from (\ref{stiffness-matrix}) we have
\bey
\V{v}^T A_{11} \V{v} &=& \sum_{K \in \mathcal{T}_h} |K| \;( \nabla v^h )^T \; \mathbb{D}_K \; \nabla v^h
+ \int_\Omega v^h \; (\V{b} \cdot \nabla v^h )d \V{x}
+ \int_\Omega c \; (v^h)^2 d \V{x}
\nn \\
&\geq& \sum_{K \in \mathcal{T}_h} |K| \;( \nabla v^h )^T \; \mathbb{D}_K \; \nabla v^h +
\int _{\Omega} ( c - \frac{1}{2} \nabla \cdot \V{b} ) ( v^h )^2 d \V{x} \ge 0.
\nn
\eey
Moreover, from the above inequality, $\V{v}^T A_{11} \V{v} = 0$ implies $v^h = $ constant, which in turn
implies $v^h = 0$ due to the fact that $v^h \in U_0$. From these, we know that $A_{11}$ is positive definite.
Since $A_{11}$ is a Z-matrix, so it is an $M$-matrix.

(4) Finally, we show matrix $A$ is an M-matrix by showing the inverse of $A$ is positive.
From (\ref{fem-matrix}), the inverse of $A$ is given by
\beq
A^{-1} = \left [\begin{array}{cc} A_{11}^{-1} & - A_{11}^{-1} A_{12} \\ 0 & I \end{array} \right ].
\eeq
Using the fact that $A_{11}^{-1} \geq 0$ and $A_{12}\le 0$, then $A^{-1} \geq 0$, which, together
with the fact that $A$ is a Z-matrix, implies that $A$ is an M-matrix.
\end{proof}

\begin{rem}
\label{rem4.1}
Loosely speaking, (\ref{thm4.1-1}) requires
\beq
\cos(\alpha_{ij,\mathbb{D}_K^{-1}}) \geq \mathcal{O}(h \|\V{b}\|_\infty)
+ \mathcal{O}(h^2 \|\V{c}\|_\infty)
\label{meshcon-1}
\eeq
or
\beq
0 < \alpha_{ij,\mathbb{D}_K^{-1}} \leq \frac{\pi}{2}-\mathcal{O}(h \|\V{b}\|_\infty)
- \mathcal{O}(h^2 \|\V{c}\|_\infty)
\label{meshcon-2}
\eeq
for all dihedral angles, where $h = \max\limits_{K \in \mathcal{T}_h} h_K$ is the maximum element size.
In other words, {\em if the mesh is $\mathcal{O}(h)$-acute in the metric $\mathbb{D}^{-1}$
for the case $\V{b} \not\equiv 0$ or $\mathcal{O}(h^2)$-acute in the metric $\mathbb{D}^{-1}$
for the case $\V{b} \equiv 0$ and $c \not\equiv 0$, then the linear finite element solution of (\ref{bvp-pde})
and (\ref{bvp-bc}) satisfies a DMP}.
\qed
\end{rem}

\begin{rem}
\label{rem4.2}
When no convection and reaction terms are involved,
condition (\ref{thm4.1-1}) reduces to the nonobtuse angle condition of \cite{CR73} and
the anisotropic nonobtuse angle condition of \cite{LH10}
for isotropic and anisotropic diffusion problems, respectively.
Moreover, the condition is consistent with the DMP conditions obtained by Brandts et al. \cite{BKK08}
and Wang and Zhang \cite{WaZh11} for isotropic diffusion-reaction or diffusion-convection-reaction
problems. Thus, (\ref{thm4.1-1}) can be viewed as a generalization of those existing results
to anisotropic diffusion-convection-reaction problems.
\qed
\end{rem}

\begin{rem}
\label{rem4.3}
Condition (\ref{thm4.1-1}) can be rewritten into a more friendly form to mesh generation.
Indeed, combining (\ref{dihedral-2}) with (\ref{thm4.1-1}) we have
\bey
&& \frac{h_i^K}{\lambda_{min}(\mathbb{D}_K)} \cdot \frac{\|\V{b}\|_{\infty, K} }{(d+1)}
+ \frac{h_i^K h_j^K}{\lambda_{min}(\mathbb{D}_K)}\cdot \frac{\|c\|_{\infty, K} }{(d+1)(d+2)}
\nn \\
&& \quad 
\le - \; \frac{\V{q}_i^T \mathbb{D}_K \V{q}_j}{\sqrt{\V{q}_i^T \mathbb{D}_K \V{q}_i} \sqrt{\V{q}_j^T \mathbb{D}_K \V{q}_j}} .
\label{rem4.3-1}
\eey
Denote the reference element and its $\V{q}$-vectors by $\hat{K}$ and $\hat{\V{q}}_i$ ($i=1, ..., d+1$),
respectively. Let $F_K$ and $F_K'$ be the affine mapping from $\hat{K}$ to $K$ and its Jacobian matrix.
As for (\ref{relation-1}), it is not difficult to show that
\[
\V{q}_i = (F_K')^{-1} \hat{\V{q}}_i, \quad i = 1, ..., d+1.
\]
Inserting this into (\ref{rem4.3-1}) leads to
\bey
&& \frac{h_i^K}{\lambda_{min}(\mathbb{D}_K)} \cdot \frac{\|\V{b}\|_{\infty, K} }{(d+1)}
+ \frac{h_i^K h_j^K}{\lambda_{min}(\mathbb{D}_K)}\cdot \frac{\|c\|_{\infty, K} }{(d+1)(d+2)}
\nn \\
&& \quad \le
- \; \frac{\hat{\V{q}}_i^T (F_K')^{-T} \mathbb{D}_K (F_K')^{-1} \hat{\V{q}}_j}
{\left (\hat{\V{q}}_i^T (F_K')^{-T} \mathbb{D}_K (F_K')^{-1} \hat{\V{q}}_i\right )^{\frac 1 2}
\left (\hat{\V{q}}_j^T (F_K')^{-T} \mathbb{D}_K (F_K')^{-1} \hat{\V{q}}_j\right )^{\frac 1 2}} .
\label{rem4.3-2}
\eey

We now consider so-called $M$-uniform meshes for a given tensor $M = M(\V{x})$ which is
assumed to be a symmetric and positive definite $d\times d$ matrix for any $\V{x} \in \overline{\Omega}
= \Omega\cup \partial \Omega$.
These meshes are approximately uniform in the metric specified by $M$. They are defined
(e.g., see \cite{Hua06,HR11}) as meshes satisfying
\beq
(F_K')^{-T} M_K^{-1} (F_K')^{-1} = \left (\frac{\sigma_h}{N}\right )^{-\frac{2}{d}} I,\quad
\forall K \in \mathcal{T}_h 
\label{M-uniform-1}
\eeq
where
\[
M_K = \frac{1}{|K|} \int_{K} M(\V{x}) d \V{x},\quad
\sigma_h = \sum_{K} |K| \; \sqrt{\text{det}(M_K)},
\]
and $\text{det}(M_K)$ is the determinant of $M_K$. $M$-uniform meshes (or more practically,
almost $M$-uniform meshes)
can be generated using a variety of techniques including blue refinement,
directional refinement, Delaunay-type triangulation, front advancing, bubble packing,
local refinement and modification, and variational mesh generation; see references in \cite{LH10}.

When the metric tensor is chosen as $M_K = \mathbb{D}_K^{-1}$ or, more generally,
\beq
M_K = \theta_K \mathbb{D}_K^{-1},
\label{M-1}
\eeq
where $\theta_K$ is a scalar, piecewise constant function, the corresponding $M$-uniform meshes will be
referred to as {\em $\mathbb{D}^{-1}$-uniform meshes}. For those meshes, (\ref{M-uniform-1}) becomes
\[
(F_K')^{-T} \mathbb{D}_K (F_K')^{-1} = \theta_K \left (\frac{\sigma_h}{N}\right )^{-\frac{2}{d}} I,\quad
\forall K \in \mathcal{T}_h .
\]
Inserting this into (\ref{rem4.3-2}, we get
\[
\frac{h_i^K}{\lambda_{min}(\mathbb{D}_K)} \cdot \frac{\|\V{b}\|_{\infty, K} }{(d+1)}
+ \frac{h_i^K h_j^K}{\lambda_{min}(\mathbb{D}_K)}\cdot \frac{\|c\|_{\infty, K} }{(d+1)(d+2)}
\le \cos(\hat{\alpha}_{ij}),
\]
where $\hat{\alpha}_{ij}$ is a dihedral angle of $\hat{K}$. If the reference element $\hat{K}$ is chosen as a regular
$d$-dimensional simplex, we have $\cos(\hat{\alpha}_{ij}) = 1/d$, and the above inequality becomes
\beq
\frac{h_i^K}{\lambda_{min}(\mathbb{D}_K)} \cdot \frac{\|\V{b}\|_{\infty, K} }{(d+1)}
+ \frac{h_i^K h_j^K}{\lambda_{min}(\mathbb{D}_K)}\cdot \frac{\|c\|_{\infty, K} }{(d+1)(d+2)}
\le \frac{1}{d} .
\label{rem4.3-3}
\eeq
This inequality holds if the maximum element size satisfies
\beq
h \|\V{b}\|_{\infty} + \frac{1}{d+2} h^2 \|c \|_{\infty} \le \frac{d+1}{d} 
\min\limits_{\V{x} \in \Omega \cup \partial \Omega} \lambda_{min}(\mathbb{D}(\V{x})) .
\label{rem4.3-4}
\eeq

It should be emphasized that (\ref{rem4.3-4}) is generally too conservative to be useful
in practical computation. However, it does show that
{\em when the reference element is chosen
as a regular $d$-dimensional simplex, a sufficiently fine (with $h$ satisfying (\ref{rem4.3-4})),
$\mathbb{D}^{-1}$-uniform mesh
satisfies the condition (\ref{thm4.1-1}). In other words, a mesh satisfying (\ref{thm4.1-1}) can be
obtained by refining a $\mathbb{D}^{-1}$-uniform mesh.}
\qed
\end{rem}

It is known that the acute or nonobtuse angle condition can be replaced by the weaker, so-called Delaunay
condition in two dimensions for a linear finite element solution to satisfy a DMP; e.g., see Strang and Fix \cite{SF73}
for the anisotropic diffusion case and Huang \cite{Hua10} for the anisotropic diffusion case.
In the current situation with convection and reaction terms, a similar weaker condition can also be obtained
in two dimensions. The argument is almost the same as that of Theorem \ref{thm4.1} except that
Step (2) of the proof needs to be fine-tuned. Let $e_{i j}$ be the edge connecting vertices $\V{a}_i$
and $\V{a}_j$ ($i = 1, ..., N_{vi},\; j = 1, ..., N, \; i \neq j$).
Denote the two elements sharing $e_{i,j}$ by $K$ and $K'$. From Step (2), we have
\begin{align}
a_{ij} & \le |K| (\nabla \phi_i|_K)^T \mathbb{D}_K \nabla \phi_j|_K
+ \frac{|K| \; \|\V{b}\|_{\infty, K}}{h_j^K (d+1)} + \frac{|K| \; \|c\|_{\infty, K}}{(d+1)(d+2)}
\nn \\
& \quad + |K'| (\nabla \phi_i|_{K'})^T \mathbb{D}_{K'} \nabla \phi_j|_{K'}
+ \frac{|K'| \; \|\V{b}\|_{\infty, K'}}{h_j^{K'} (d+1)} + \frac{|K'| \; \|c\|_{\infty, K'}}{(d+1)(d+2)}
\nn \\
& = - \frac{\mbox{det}(\mathbb{D}_K)^{\frac 1 2}}{2} \cot(\alpha_{ij, \mathbb{D}_K^{-1}})
- \frac{\mbox{det}(\mathbb{D}_{K'})^{\frac 1 2}}{2} \cot(\alpha_{ij, \mathbb{D}_{K'}^{-1}})
\nn \\
& \quad + \frac{|K| \; \|\V{b}\|_{\infty, K}}{h_j^K (d+1)} + \frac{|K| \; \|c\|_{\infty, K}}{(d+1)(d+2)}
+ \frac{|K'| \; \|\V{b}\|_{\infty, K'}}{h_j^{K'} (d+1)} + \frac{|K'| \; \|c\|_{\infty, K'}}{(d+1)(d+2)} ,
\tag{Lemma \ref{lem3.3}}
\end{align}
where $\alpha_{ij, \mathbb{D}_K^{-1}}$ and $\alpha_{ij, \mathbb{D}_{K'}^{-1}}$ are the angles of $K$ and $K'$,
respectively, that face the common edge $e_{i j}$.
From this we can conclude that the linear finite element solution in 2D satisfies a DMP if the mesh satisfies
\bey
& & \frac{|K| \; \|\V{b}\|_{\infty, K}}{h_j^K (d+1)} + \frac{|K| \; \|c\|_{\infty, K}}{(d+1)(d+2)}
+ \frac{|K'| \; \|\V{b}\|_{\infty, K'}}{h_j^{K'} (d+1)} + \frac{|K'| \; \|c\|_{\infty, K'}}{(d+1)(d+2)}
\nn \\
& \le & \frac{\mbox{det}(\mathbb{D}_K)^{\frac 1 2}}{2} \cot(\alpha_{ij, \mathbb{D}_K^{-1}})
+ \frac{\mbox{det}(\mathbb{D}_{K'})^{\frac 1 2}}{2} \cot(\alpha_{ij, \mathbb{D}_{K'}^{-1}})
\label{meshcon-3}
\eey
for all internal edges. Following \cite{Hua10}, we can rewrite the above inequality as
\bey
&& 0 < \frac{1}{2} \left [ \frac{}{} \alpha_{ij, \mathbb{D}_K^{-1}} + \alpha_{ij, \mathbb{D}_{K'}^{-1}}
\right . \nn \\
&& \quad
+ \mbox{ arccot}\left (\sqrt{\frac{\mbox{ det}(\mathbb{D}_{K'})}{\mbox{ det}(\mathbb{D}_{K})}}
\cot (\alpha_{ij, \mathbb{D}_{K'}^{-1}} ) - \frac{2\; C(K,K',j) }{\sqrt{\mbox{ det}(\mathbb{D}_{K})}} \right )
\nn \\
&& \quad \left . + \mbox{ arccot}\left (\sqrt{\frac{\mbox{ det}(\mathbb{D}_{K})}{\mbox{ det}(\mathbb{D}_{K'})}}
\cot (\alpha_{ij, \mathbb{D}_{K}^{-1}} ) - \frac{2\; C(K,K',j) }{\sqrt{\mbox{ det}(\mathbb{D}_{K'})}} \right )
\frac{}{} \right ] \le \pi,
\label{meshcon-4}
\eey
where
\bey
C(K,K',j) & & = \frac{|K| \; \|\V{b}\|_{\infty, K}}{h_j^K (d+1)} + \frac{|K| \; \|c\|_{\infty, K}}{(d+1)(d+2)}
+ \frac{|K'| \; \|\V{b}\|_{\infty, K'}}{h_j^{K'} (d+1)}
\nn \\
&& \quad + \; \frac{|K'| \; \|c\|_{\infty, K'}}{(d+1)(d+2)} .
\label{meshcon-5}
\eey

The following theorem summarizes the above analysis.

\begin{thm}
\label{thm4.2}
If (\ref{meshcon-4}) holds for all internal edges of the simplicial mesh $\mathcal{T}_h$,
then the linear finite element scheme (\ref{disc-1}) for BVP (\ref{bvp-pde}) and (\ref{bvp-bc}) in two dimensions
satisfies DMP.
\end{thm}

Loosely speaking, (\ref{meshcon-4}) can be written as
\bey
&& 0 < \frac{1}{2} \left [ \frac{}{} \alpha_{ij, \mathbb{D}_K^{-1}} + \alpha_{ij, \mathbb{D}_{K'}^{-1}}
+ \mbox{ arccot}\left (\sqrt{\frac{\mbox{ det}(\mathbb{D}_{K'})}{\mbox{ det}(\mathbb{D}_{K})}}
\cot (\alpha_{ij, \mathbb{D}_{K'}^{-1}} ) \right )
\right .
\nn \\
&& \qquad \qquad \quad \left . + \mbox{ arccot}\left (\sqrt{\frac{\mbox{ det}(\mathbb{D}_{K})}{\mbox{ det}(\mathbb{D}_{K'})}}
\cot (\alpha_{ij, \mathbb{D}_{K}^{-1}} ) \right )
\frac{}{} \right ] 
\nn \\
&& \quad \le \pi - \mathcal{O}(h \|\V{b}\|_\infty) - \mathcal{O}(h^2 \|c\|_\infty) .
\label{meshcon-6}
\eey

\begin{rem}
\label{rem4.4}
For the case where $\mathbb{D} = I$, $\V{b} = 0$, and $c = 0$, it is easy to see that
(\ref{meshcon-6}) reduces to the Delaunay condition: $\alpha_{ij,K} + \alpha_{ij,K'} \le \pi$.
Moreover, for the case without convection and reaction terms, (\ref{meshcon-6}) gives
the Delaunay-type mesh condition obtained by Huang \cite{Hua10} for two dimensional
anisotropic diffusion problems.
\qed
\end{rem}

\section{Numerical examples}
In this section we present numerical results obtained for four 2D examples to verify the
mesh condition (\ref{thm4.1-1}) and (\ref{meshcon-6}). In all but Example \ref{exam5.4},
the convection vector $\V{b}$ is taken as a constant vector with equal, positive $x$ and $y$ components,
i.e., $\V{b} = \|\V{b}\|_\infty (1, 1)^T$.

\begin{figure}[thb]
\centering
\includegraphics[width=3in]{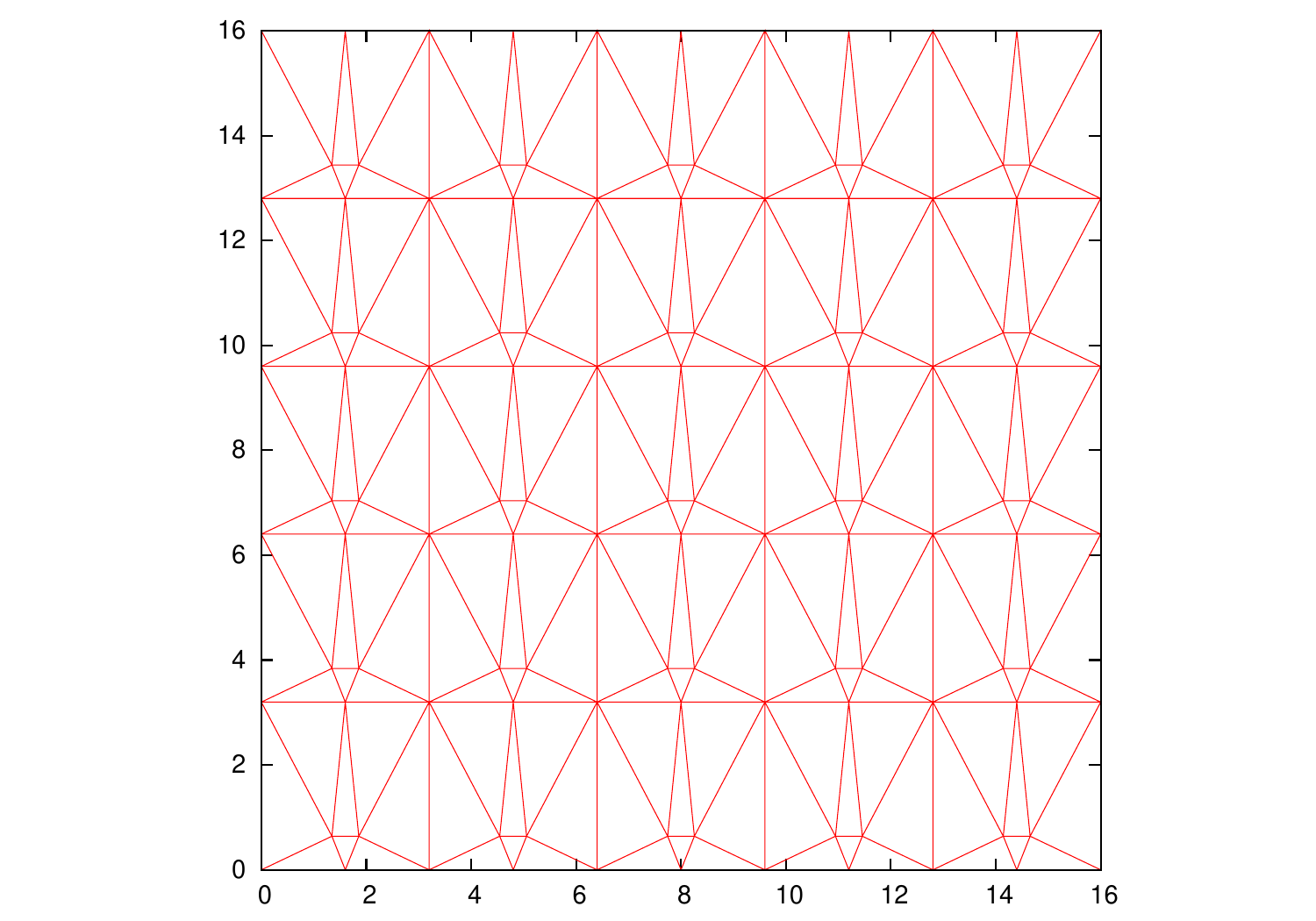}
\caption{A typical mesh ($N=200$) used for Example \ref{exam5.1}.}
\label{Exa5.1-1}
\end{figure}

\begin{figure}[thb]
\centering
\hbox{
\hspace{17mm}
\begin{minipage}[b]{2.5in}
\centerline{(a): $N=9800$}
\includegraphics[width=2.5in]{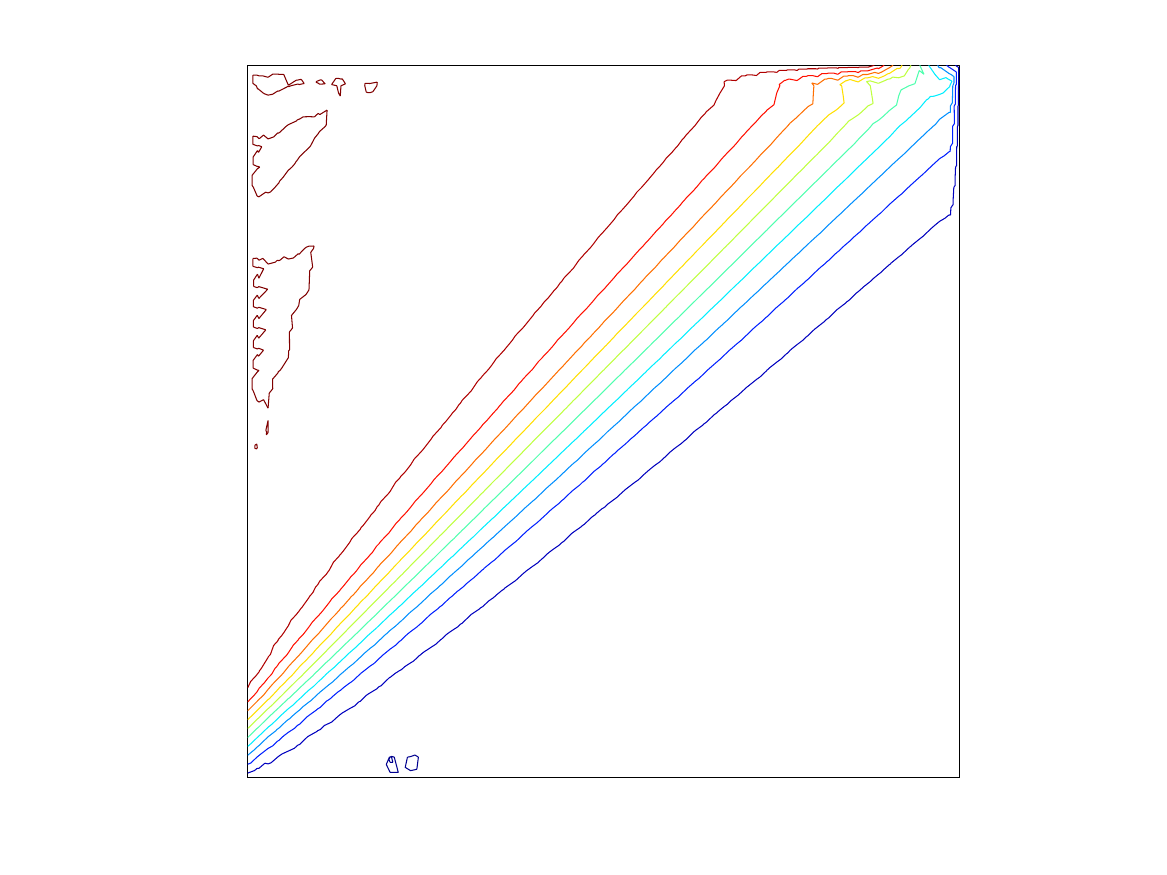}
\end{minipage}
\hspace{-10mm}
\begin{minipage}[b]{2.5in}
\centerline{(b): $N=20000$}
\includegraphics[width=2.5in]{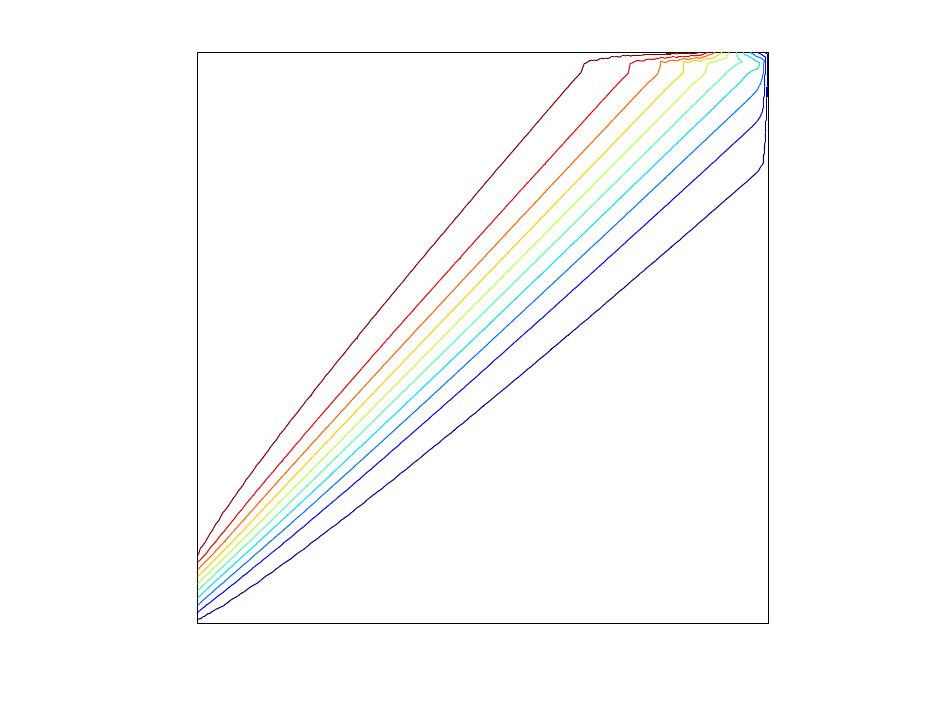}
\end{minipage}
}
\caption{Contours of the linear finite element solutions for Example \ref{exam5.1}.}
\label{Exa5.1-2}
\end{figure}

\begin{figure}[thb]
\centering
\hbox{
\hspace{8mm}
\begin{minipage}[b]{3in}
\centerline{(a): $\V{b} = (10,10)^T$}
\includegraphics[width=3in]{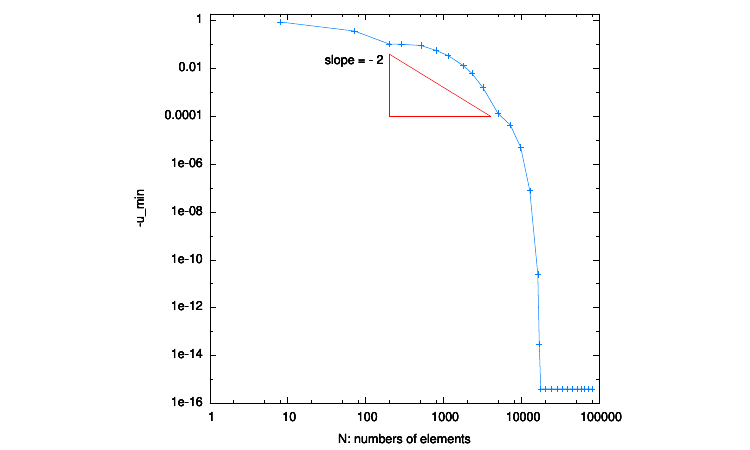}
\end{minipage}
\hspace{-20mm}
\begin{minipage}[b]{3in}
\centerline{(b): $N = 3200$}
\includegraphics[width=3in]{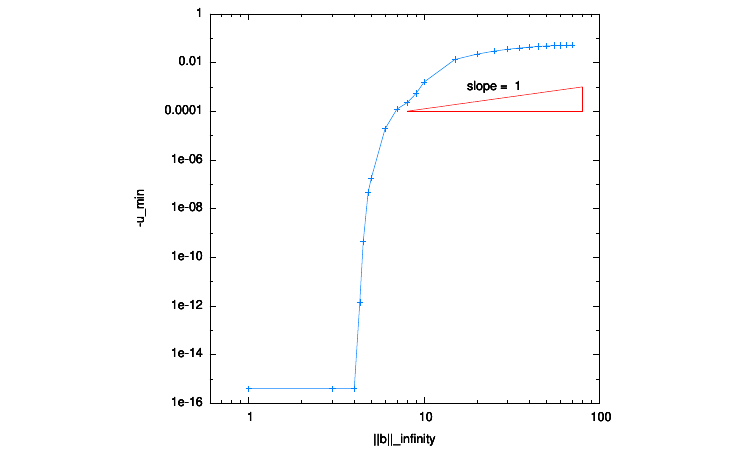}
\end{minipage}
}
\caption{The undershoot, $-u_{min}$, is plotted as  a function of the number of elements $N$ in (a) 
and as a function of $||\V{b}||_{\infty}$ in (b) for Example \ref{exam5.1}.}
\label{Exa5.1-3}
\end{figure}

\begin{exam}
\label{exam5.1}
The first example is in the form of (\ref{bvp-pde}) and (\ref{bvp-bc}),
and the coefficients are given as
\begin{align*}
& c = 0, \;\; f = 0,\;\; g(x,0)=g(16,y) = 0,
\\
& g(0,y)= \begin{cases} 0.5y, & \text{ for } 0 \le y < 2 \\
1, & \text{ for } 2 \le y \le 16 \end{cases} 
\\
& g(x,16)= \begin{cases} 1, & \text{ for } 0 \le x \le 14 \\
8-0.5x, & \text{ for } 14 < x \le 16. \end{cases}
\end{align*}
For this example, the diffusion matrix is taken as the identity matrix, i.e., $\mathbb{D}=I$.
This is an isotropic homogeneous diffusion problem.
Note that the example satisfies the maximum principle and its solution stays between $0$ and $1$. 

An acute-type mesh is used in the computation.
Such a mesh is obtained by partitioning  each square element of a uniform mesh into eight triangles with acute angles;
see Fig.~\ref{Exa5.1-1}. The maximum angle of the mesh is $0.49\pi$ and thus
condition (\ref{thm4.1-1}) holds when the mesh size is sufficiently small.

Fig.~\ref{Exa5.1-2} shows the contours of the linear finite element solutions obtained for
$N=9800$ and $N=20000$. ($N$ is the number of elements.)
There are no undershoot nor overshoot for $N=20000$ whereas both undershoots and overshoots occur
for the case with $N=9800$.
In Fig.~\ref{Exa5.1-3}(a), $-u_{min}$ is shown as functions of the number of elements $N$.
From the figure one can see that $-u_{min}$ decreases as the mesh is refined and
the decrease rate is about quadratic initially and then exponential near $N=10000$.
Moreover, $-u_{min}$ becomes zero (more precisely, at the level of roundoff error)
after around $N = 17000$. This is consistent with Theorem \ref{thm4.1}
which states that there are no undershoot nor overshoot when the mesh size is sufficiently small.

To further verify Theorem \ref{thm4.1}, we fix the number of elements at $N=3200$ and let
$||\V{b}||_{\infty}$ vary. Quantity $-u_{min}$ is plotted in Fig.~\ref{Exa5.1-3}(b) as a function of
$||\V{b}||_{\infty}$. From the figure, we can see that there is no undershoot until $||\V{b}||_{\infty} \approx 4$.
Then $-u_{min}$ increases exponentially until $||\V{b}||_{\infty} \approx 20$ where
the increase rate is about linear as $||\V{b}||_{\infty}$ increases.

Finally, it is pointed out that a similar behavior can be observed for the overshoot. The results are omitted
here to save space.
\qed
\end{exam}


\begin{figure}[thb]
\centering
\includegraphics[width=3in]{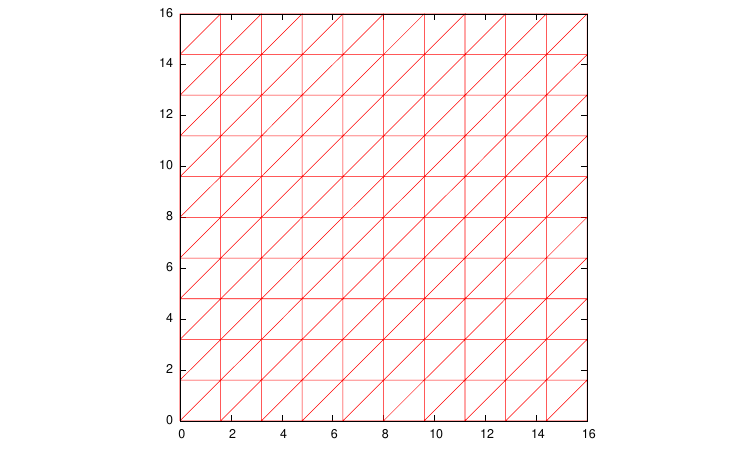}
\caption{A typical mesh ($N = 200$) used for Example \ref{exam5.2} (with anisotropic $\mathbb{D}$). }
\label{Exa5.2-1}
\end{figure}

\begin{figure}[thb]
\centering
\hbox{
\hspace{17mm}
\begin{minipage}[b]{2.5in}
\centerline{(a): $N = 450$}
\includegraphics[width=2.5in]{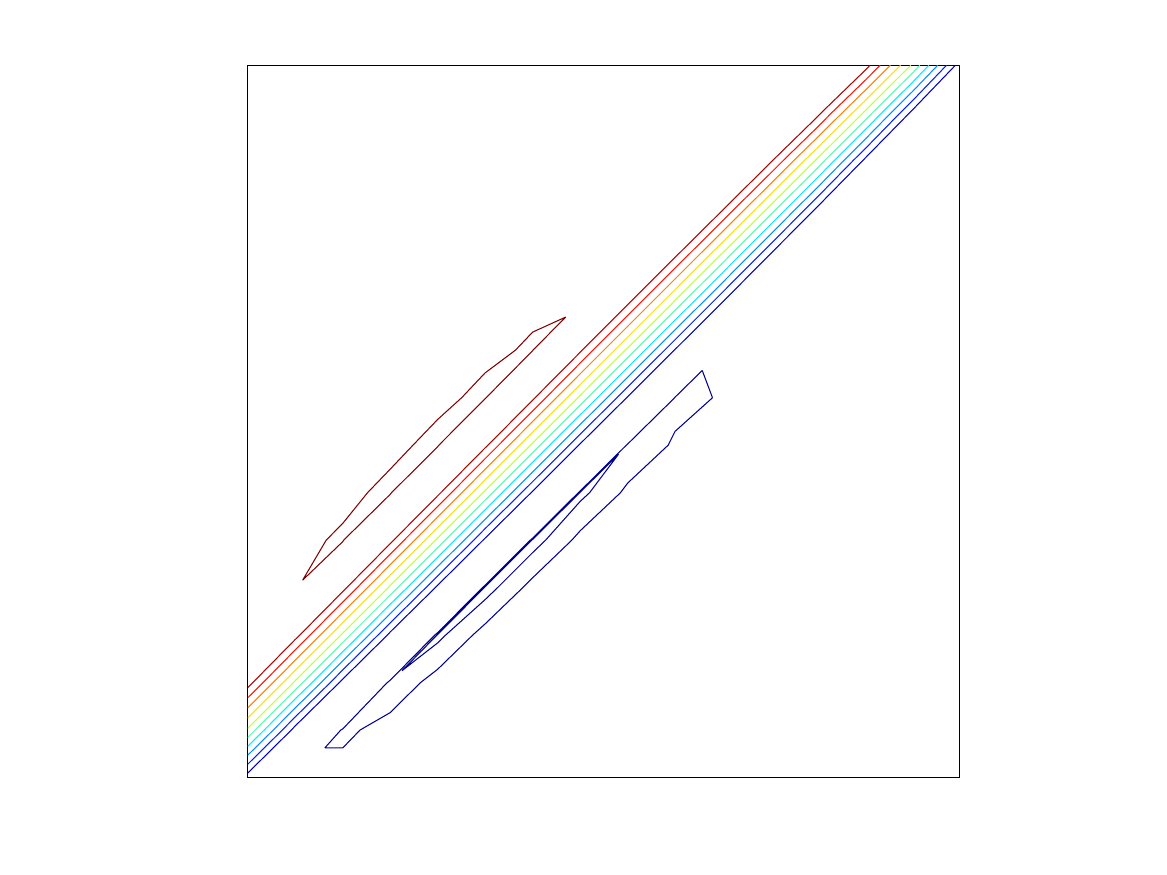}
\end{minipage}
\hspace{-10mm}
\begin{minipage}[b]{2.5in}
\centerline{(b): $N = 7200$}
\includegraphics[width=2.5in]{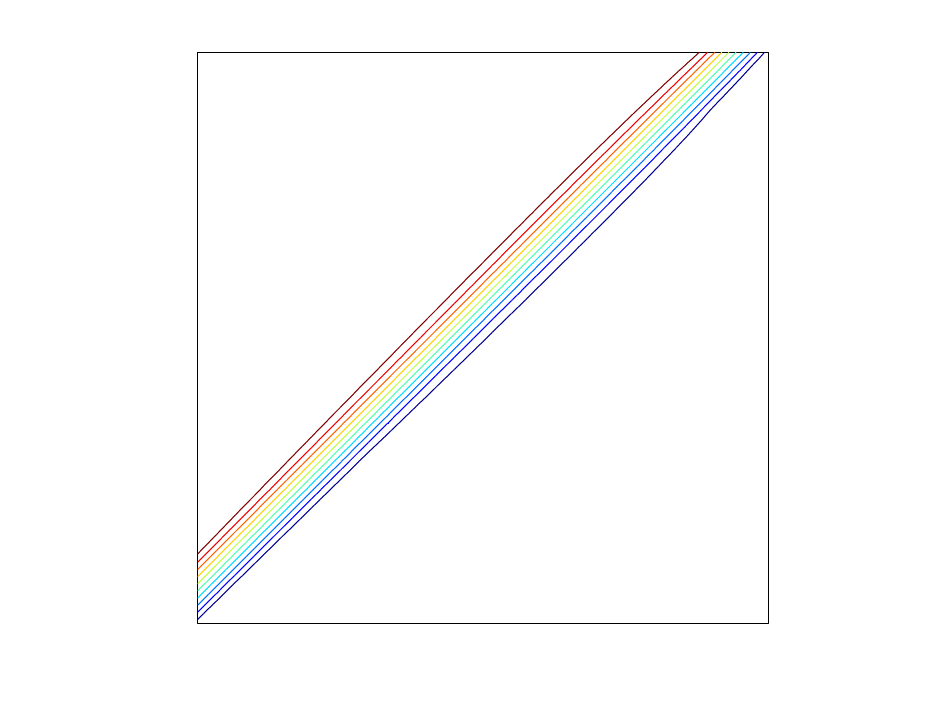}
\end{minipage}
}
\caption{Contours of linear finite element solutions for Example \ref{exam5.2}. }
\label{Exa5.2-2}
\end{figure}

\begin{figure}[thb]
\centering
\hbox{
\hspace{8mm}
\begin{minipage}[b]{3in}
\centerline{(a): $\V{b}=[200,200]^T$}
\includegraphics[width=3in]{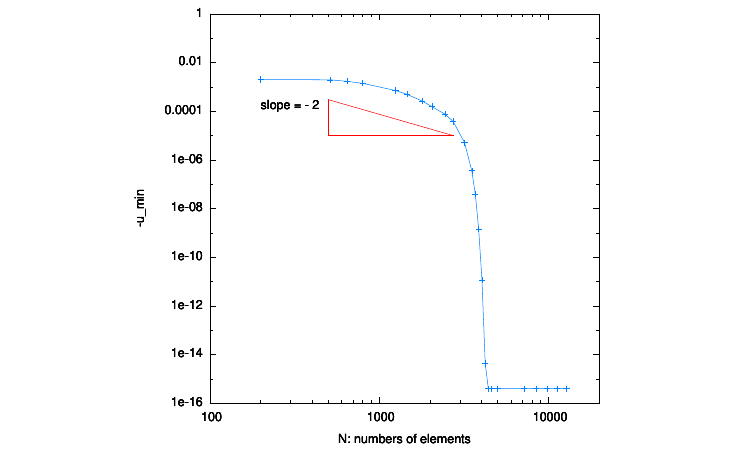}
\end{minipage}
\hspace{-20mm}
\begin{minipage}[b]{3in}
\centerline{(b): $N=3200$}
\includegraphics[width=3in]{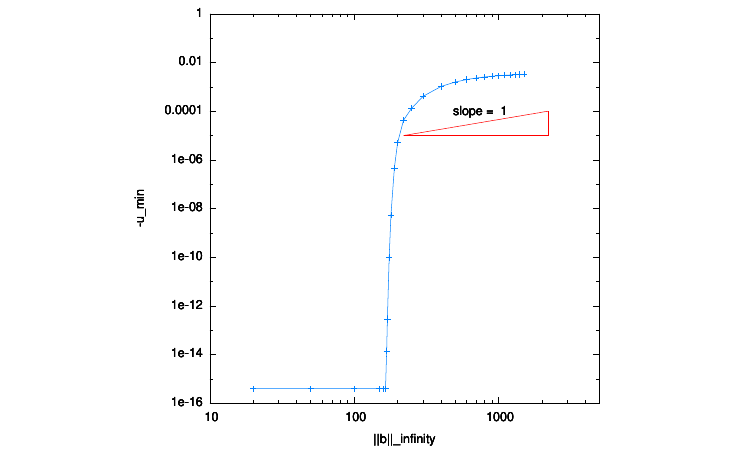}
\end{minipage}
}
\caption{The undershoot, $-u_{min}$, is plotted as  a function of the number of elements $N$ in (a) 
and as a function of $||\V{b}||_{\infty}$ in (b) for Example \ref{exam5.2}.}
\label{Exa5.2-3}
\end{figure}

\begin{exam}
\label{exam5.2}
In the second example, BVP (\ref{bvp-pde}) and (\ref{bvp-bc}) with all the coefficients being the same with 
Example \ref{exam5.1} except the diffusion matrix is used.
The diffusion matrix is taken as
\beq
\nn
\mathbb{D}(x,y) = \left ( \begin{array}{cc} 500.5 & 499.5 \\ 499.5 & 500.5 \end{array} \right ) .
\eeq
This matrix represents a homogeneous but highly anisotropic diffusion process.
A mesh with right triangle elements (see Fig.~\ref{Exa5.2-1}) is used for this example. 
Such a mesh is obtained by dividing each square element of a uniform mesh into two right triangular elements.
Although each element of the mesh is a right triangle (in the Euclidean sense), the maximum angle is $0.49\pi$
when measured in metric $\mathbb{D}^{-1}$. Thus, the mesh is of acute-type in the metric and
condition (\ref{thm4.1-1})  can be satisfied if the mesh size is sufficiently small.

Contours of linear finite element solutions are shown in Fig.~\ref{Exa5.2-2} while the undershoot
is plotted as functions of $N$ and $\|\V{b}\|_\infty$ in Fig.~\ref{Exa5.2-3}. From these results we can
observe a similar behavior of the undershoot and overshoot as in Example \ref{exam5.1}, i.e.,
they occur only for relatively coarse meshes or relatively large $\|\V{b}\|_\infty$. The behavior is
consistent with Theorem \ref{thm4.1}.
\qed
\end{exam}

\begin{figure}[thb]
\centering
\hbox{
\hspace{17mm}
\begin{minipage}[b]{2.5in}
\centerline{(a): $N = 2312$}
\includegraphics[width=2.5in]{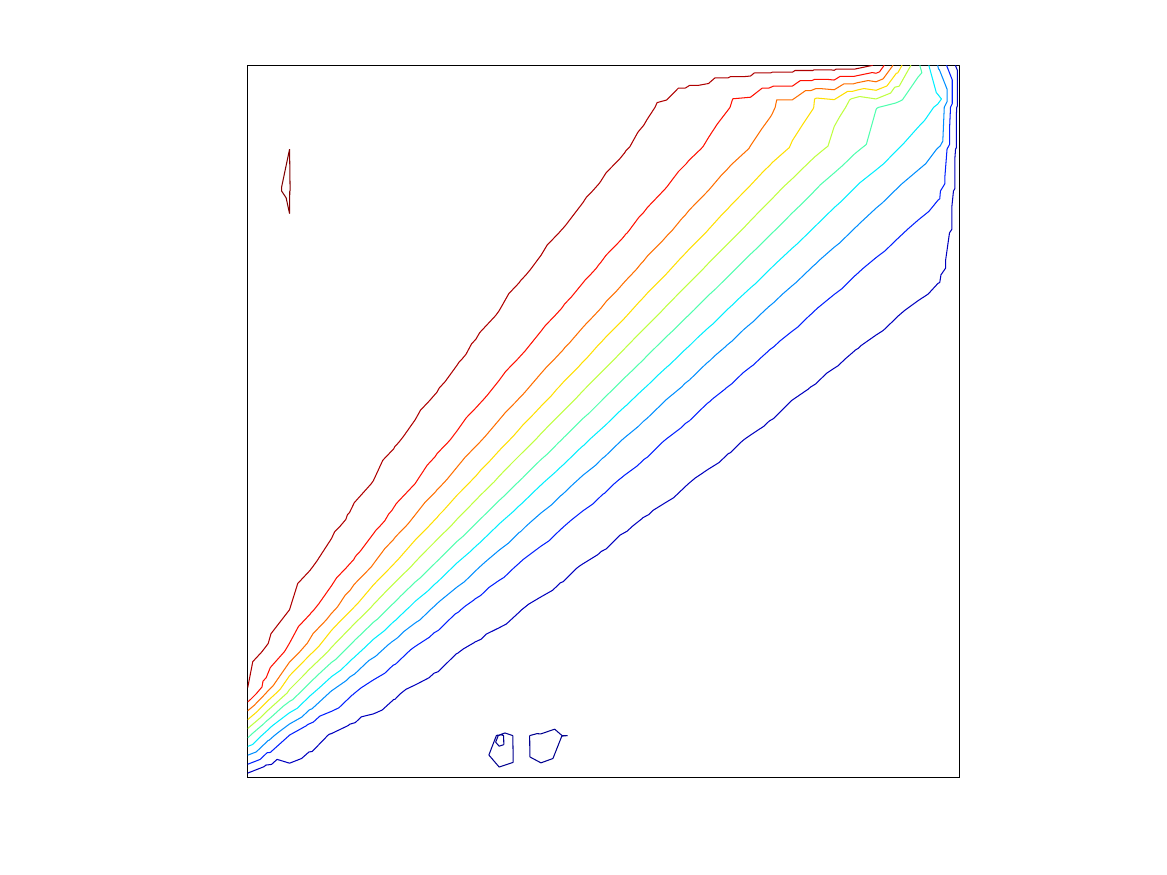}
\end{minipage}
\hspace{-10mm}
\begin{minipage}[b]{2.5in}
\centerline{(b): $N = 4232$}
\includegraphics[width=2.5in]{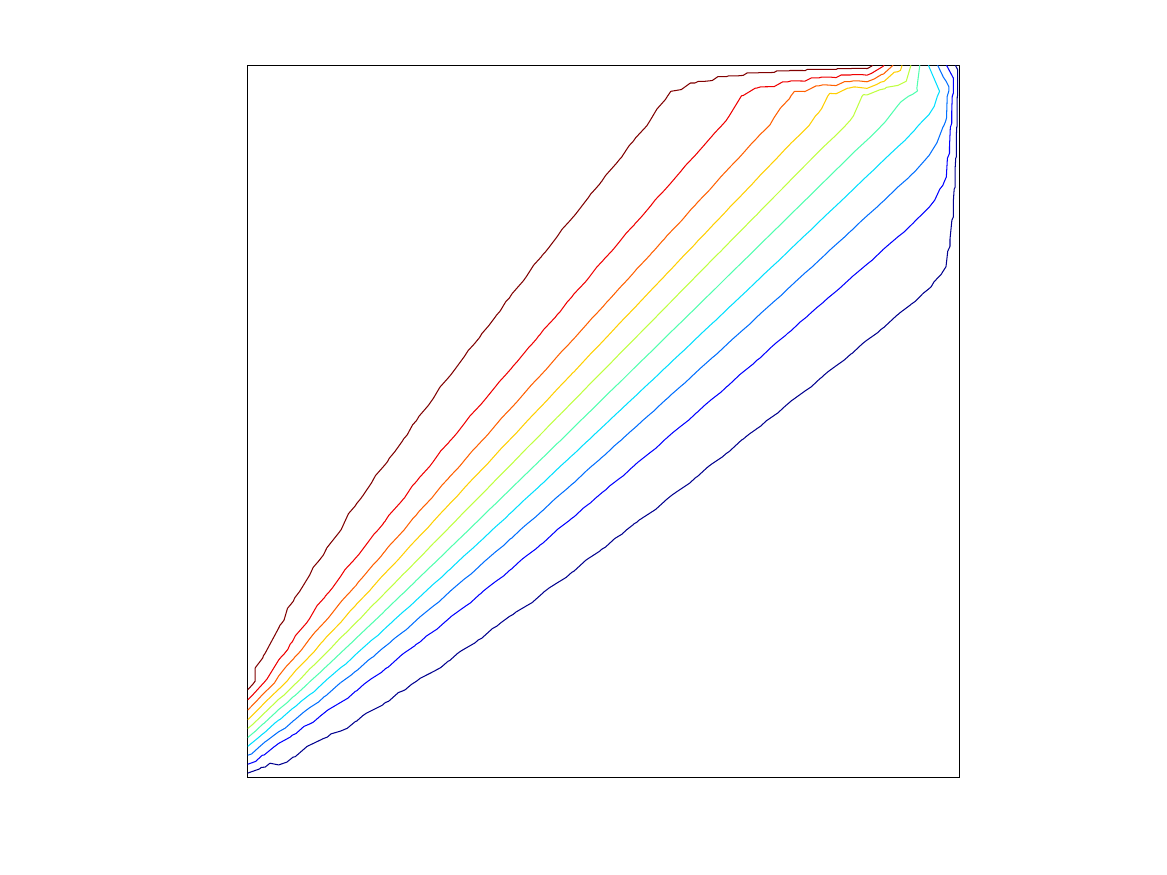}
\end{minipage}
}
\caption{Contours of linear finite element solutions for Example \ref{exam5.3}. }
\label{Exa5.3-2}
\end{figure}

\begin{figure}[thb]
\centering
\hbox{
\hspace{8mm}
\begin{minipage}[b]{3in}
\centerline{(a): $\V{b}=[200,200]^T$}
\includegraphics[width=3in]{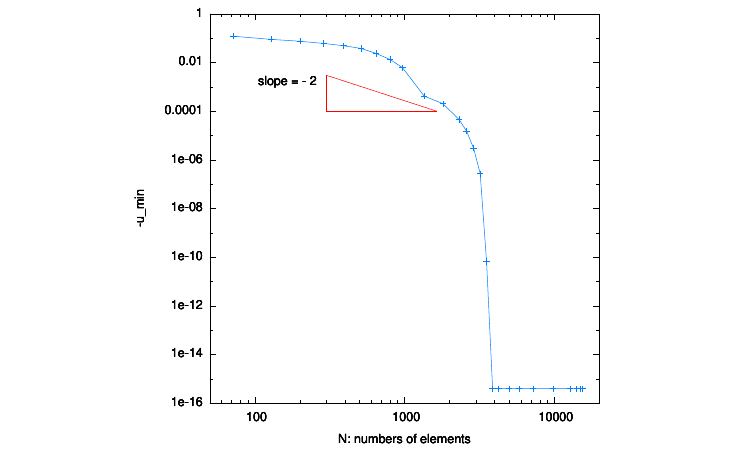}
\end{minipage}
\hspace{-20mm}
\begin{minipage}[b]{3in}
\centerline{(b): $N=3200$}
\includegraphics[width=3in]{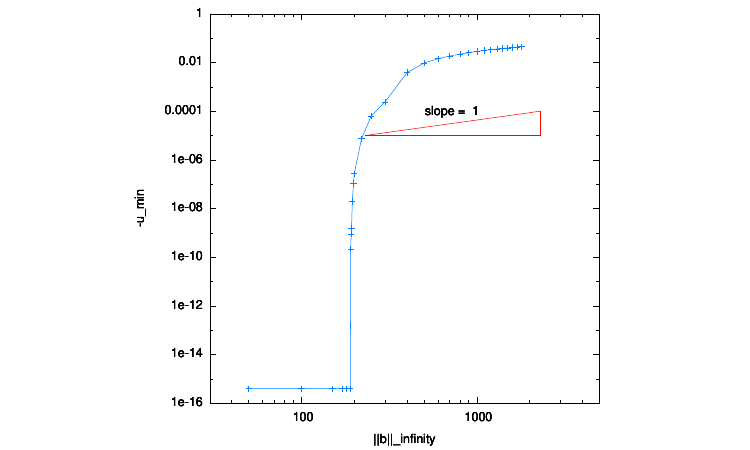}
\end{minipage}
}
\caption{The undershoot, $-u_{min}$, is plotted as  a function of the number of elements $N$ in (a) 
and as a function of $||\V{b}||_{\infty}$ in (b) for Example \ref{exam5.3}.}
\label{Exa5.3-3}
\end{figure}

\begin{exam}
\label{exam5.3}
In this example, the same BVP (\ref{bvp-pde}) and (\ref{bvp-bc}) with  Example \ref{exam5.1}
 except the diffusion matrix is used for this example.
The diffusion matrix is taken as
\beq
\nn
\mathbb{D}(x,y) = \left ( \begin{array}{cc} 50 & 12 \\ 12 & 50 \end{array} \right ) .
\eeq
This diffusion matrix has a weaker anisotropy than that in the previous example.

The mesh used for Example \ref{exam5.1} (see Fig.~\ref{Exa5.1-1}) is also used for this example. 
Recall that the mesh is acute in the Euclidean sense. 
When measured in metric $\mathbb{D}^{-1}$, however, the maximum angle of the mesh is $0.55\pi$
and the maximum sum of any pair of angles opposite a common edge is $0.97\pi$.
Thus, the mesh will satisfy (\ref{meshcon-6}) but not (\ref{thm4.1-1}) when its size is sufficiently small.

Contours of numerical solutions are shown in Fig.~\ref{Exa5.3-2} while the undershoot
is plotted as functions of $N$ and $\|\V{b}\|_\infty$ in Fig.~\ref{Exa5.3-3}. 
A similar behavior of the undershoot and overshoot can be observed as for the two previous examples.

In particular, for $N\ge 4000$ there is no undershoot or overshoot in the numerical solutions for the case
with $\V{b} =  [200, 200]^T$. This example shows that condition (\ref{meshcon-6}) is weaker
than condition (\ref{thm4.1-1}).
\qed
\end{exam}

\begin{figure}[thb]
\centering
\hbox{
\hspace{8mm}
\begin{minipage}[b]{3in}
\centerline{(a): $N = 2050$}
\includegraphics[width=3in]{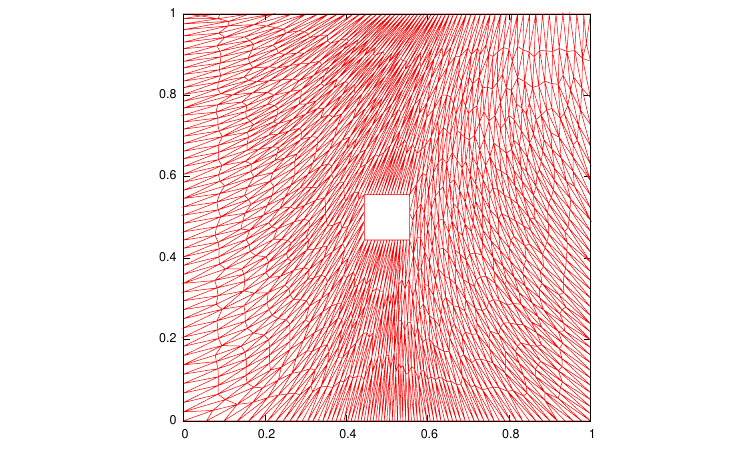}
\end{minipage}
\hspace{-20mm}
\begin{minipage}[b]{3in}
\centerline{(b): $N = 2063$}
\includegraphics[width=3in]{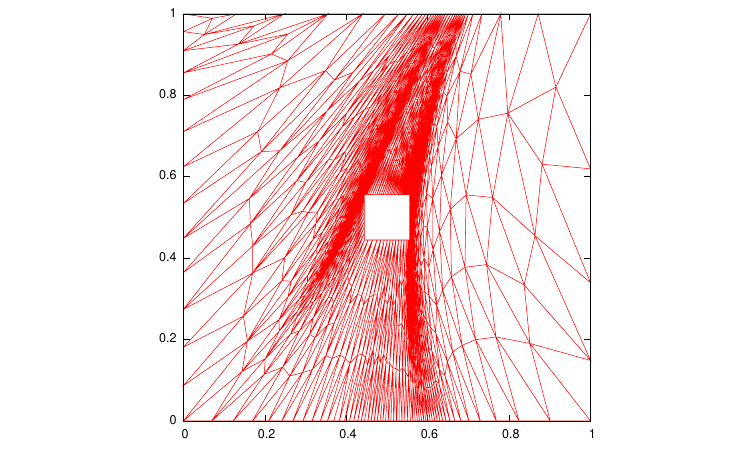}
\end{minipage}
}
\caption{Two typical (a) $M_{DMP}$ and (b) $M_{DMP+adap}$ meshes obtained for Example \ref{exam5.4}. }
\label{Exa5.4-1}
\end{figure}

\begin{figure}[thb]
\centering
\hbox{
\hspace{8mm}
\begin{minipage}[b]{3in}
\centerline{(a): $M_{DMP}$ mesh}
\includegraphics[width=3in]{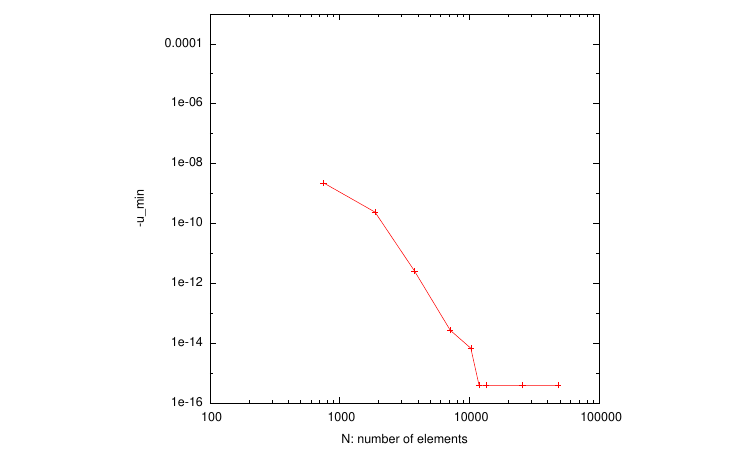}
\end{minipage}
\hspace{-20mm}
\begin{minipage}[b]{3in}
\centerline{(b): $M_{DMP+adap}$ mesh}
\includegraphics[width=3in]{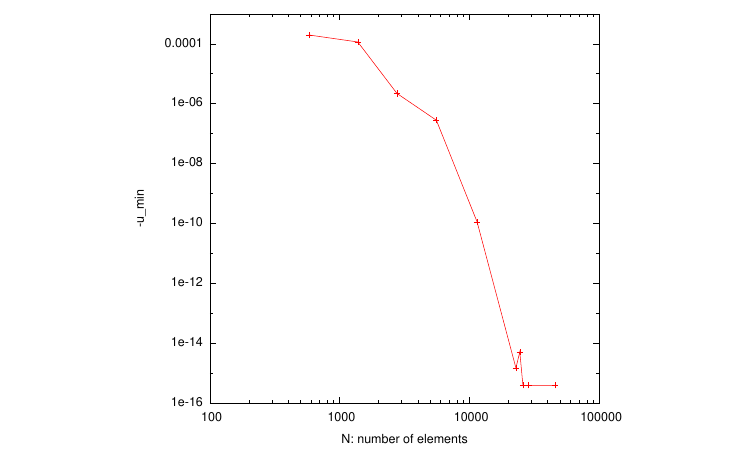}
\end{minipage}
}
\caption{The undershoot, $-u_{min}$, is plotted as  a function of the number of elements $N$ for
(a) $M_{DMP}$ and (b) $M_{DMP+adap}$ meshes for Example \ref{exam5.4}.}
\label{Exa5.4-4}
\end{figure}

\begin{figure}[thb]
\centering
\includegraphics[width=3in]{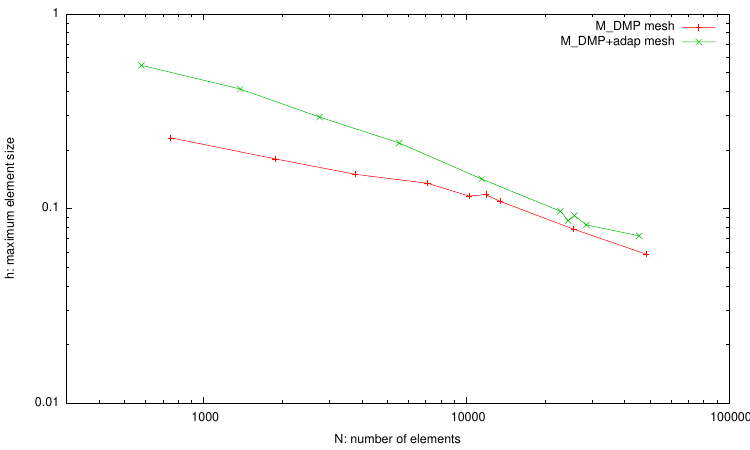}
\caption{The maximum element size is plotted as a function of the number of elements for both $M_{DMP}$ and
$M_{DMP+adap}$ meshes for Example \ref{exam5.4}.}
\label{Exa5.4-5}
\end{figure}

\begin{figure}[thb]
\centering
\hbox{
\hspace{33mm}
\begin{minipage}[b]{1.5in}
\centerline{(a): $N = 5460$}
\includegraphics[width=1.5in]{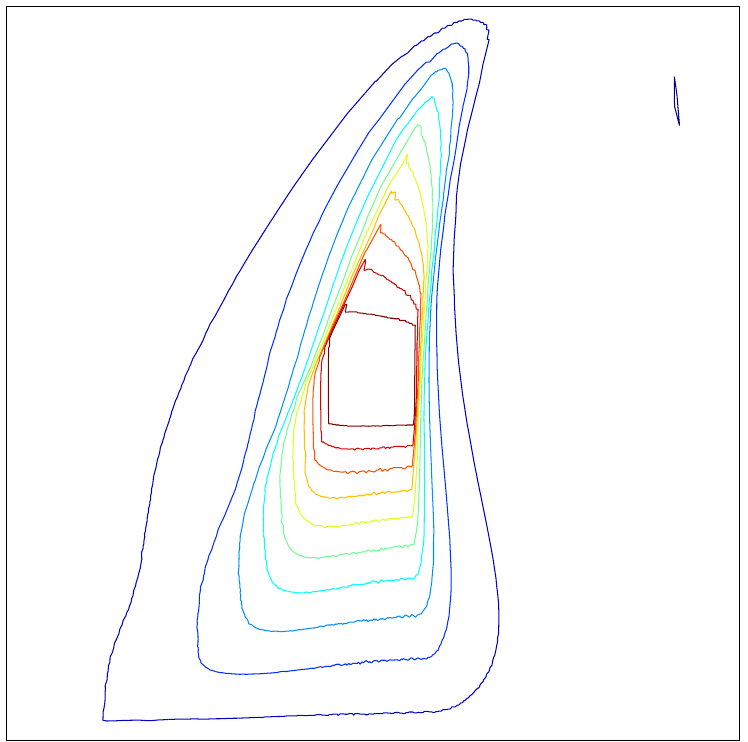}
\end{minipage}
\hspace{15mm}
\begin{minipage}[b]{1.5in}
\centerline{(b): $N = 31295$}
\includegraphics[width=1.5in]{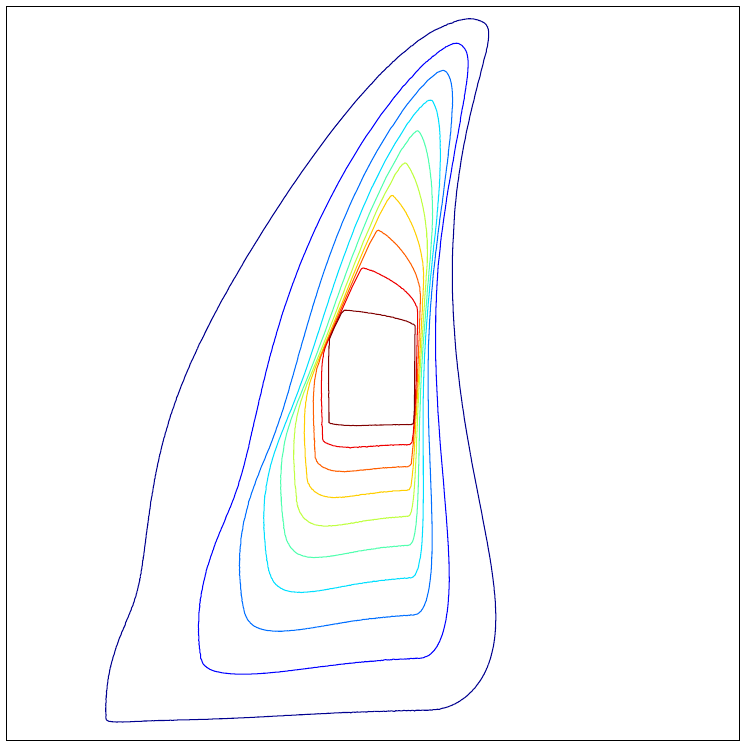}
\end{minipage}
}
\caption{Contours of linear finite element solutions obtained with $M_{DMP}$ meshes
for Example \ref{exam5.4}. }
\label{Exa5.4-2}
\end{figure}

\begin{figure}[thb]
\centering
\hbox{
\hspace{33mm}
\begin{minipage}[b]{1.5in}
\centerline{(a): $N = 5469$}
\includegraphics[width=1.5in]{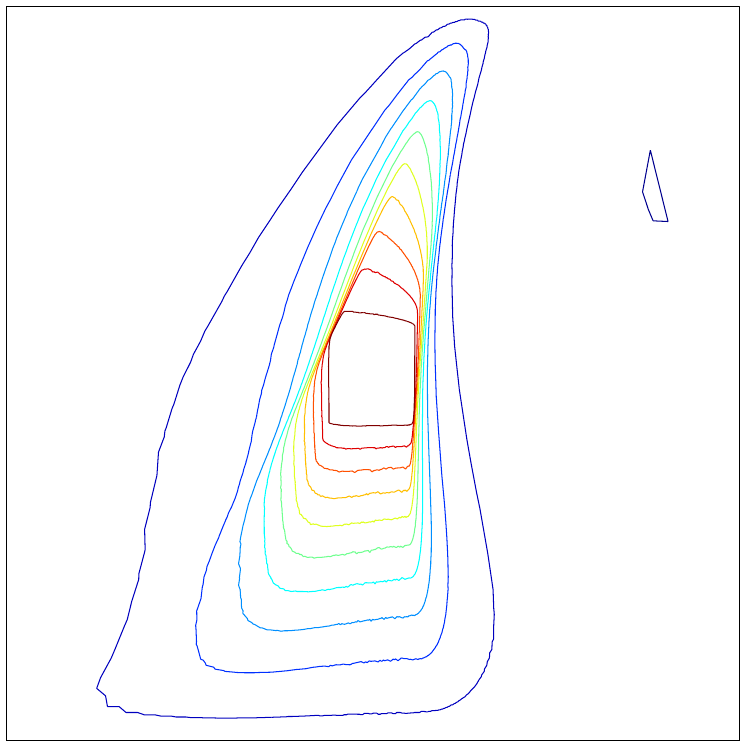}
\end{minipage}
\hspace{15mm}
\begin{minipage}[b]{1.5in}
\centerline{(b): $N = 31454$}
\includegraphics[width=1.5in]{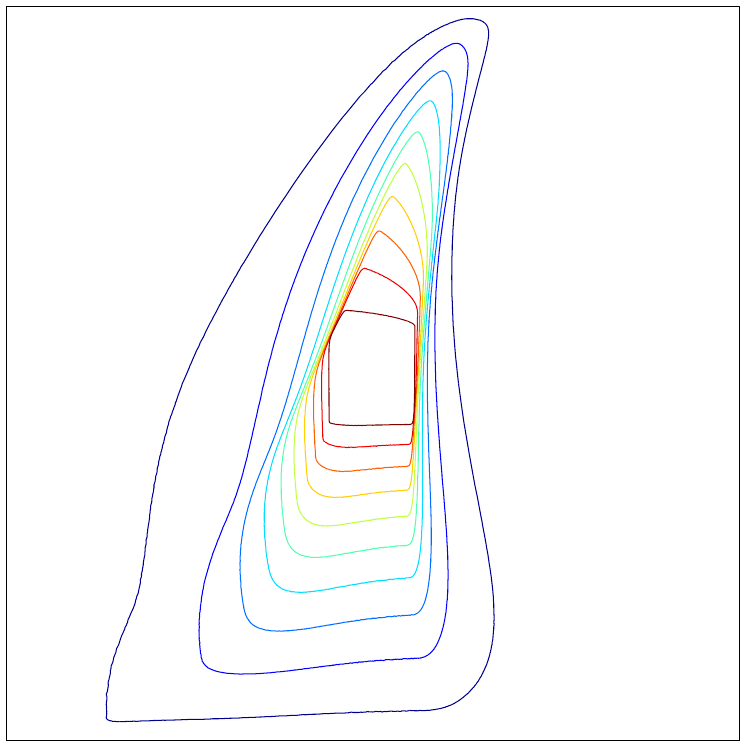}
\end{minipage}
}
\caption{Contours of linear finite element solutions obtained with $M_{DMP+adap}$ meshes
for Example \ref{exam5.4}. }
\label{Exa5.4-3}
\end{figure}

\begin{exam}
\label{exam5.4}
In the last example, we consider BVP  (\ref{bvp-pde}) and (\ref{bvp-bc}) on domain $\Omega = [0,1]^2 \backslash [\frac{4}{9},\frac{5}{9}]^2$ with 
\beq
\begin{cases}
& \V{b} = \Big [5000 (0.5-y), 5000(x-0.5) \Big]^T, \;\;  c = 100, \;\; f = 0,\\
& g = 0 \;\; \text{on} \;\; \partial \Omega_{\text{out}}, \;\; g = 2 \;\; \text{on} \;\; \partial \Omega_{\text{in}},
\end{cases}
\label{exam5.4-1}
\eeq
where $\partial \Omega_{\text{out}}$ and $\partial \Omega_{\text{in}}$ are the outer and inner boundaries of $\Omega$, respectively.
The diffusion matrix is taken as
\beq
\nn
\mathbb{D}(x,y) = 
\left ( \begin{array}{cc} \cos \alpha & -\sin \alpha \\ \sin \alpha & \cos \alpha \end{array} \right )  
\left ( \begin{array}{cc} 1000 & 0 \\ 0 & 1 \end{array} \right )  
\left ( \begin{array}{cc} \cos \alpha & \sin \alpha \\ -\sin \alpha & \cos \alpha \end{array} \right ) ,
\eeq
where $\alpha = \pi \sin(x)\cos(y)$. 
This diffusion matrix is anisotropic and heterogeneous. 
The BVP satisfies the maximum principle and its solution stays between $0$ and $2$.

We use two sets of $\mathbb{D}^{-1}$-uniform meshes (cf. Remark~\ref{rem4.3}) for this example.
The first set (referred to as $M_{DMP}$ meshes)
is $M$-uniform meshes with $M$ defined in (\ref{M-1}) and $\theta_K = 1$.
Notice that this set of meshes is completely determined by the diffusion matrix $\mathbb{D}$.
The other set of meshes (referred to as $M_{DMP+adap}$ meshes)
is $M$-uniform meshes with $M$ defined in the form (\ref{M-1}) and
$\theta_K$ determined by minimizing an interpolation error estimate. $M$ has the expression
\cite{LH10}[equation (55)] as
\beq
M_K = \left ( 1 + \frac{1}{\alpha_h}B_K \right )^{\frac{1}{2}} \sqrt{\text{det}(\mathbb{D}_K)} \;
\mathbb{D}_K^{-1},
\label{M-2}
\eeq
where 
\bey
&& B_K = \text{det}(\mathbb{D}_K)^{-\frac 1 2} \| \mathbb{D}_K^{-1} || \cdot \frac{1}{|K|}
\int_K \|\, \mathbb{D}_K |H(u)| \,\|^2 d \V{x},
\nn \\
&& \alpha_h = \left (\frac{1}{|\Omega|} \sum\limits_K |K| \sqrt{B_K} \right )^2,
\nn
\eey
$\|\cdot \|$ is the $l_2$ matrix norm, $|H(u)| = \sqrt{H(u)^2}$, and $H(u)$ denotes
the Hessian of the exact solution $u$. In the computation, the integral is calculated with
a Gaussian quadrature rule and the Hessian is replaced
by approximations obtained with a Hessian recovery technique \cite{LH10}
which employs piecewise quadratic polynomials fitting in least-squares
sense to nodal values of the currently available computed solution.

An iterative procedure is used for solving this example. It involves three basic steps, solving the BVP
on the current mesh, computing the metric tensor, and generating a new mesh.
In our computation, each run is stopped after 10 iterations.
We have found that there is very little improvement in the computed solution after 10 iterations.
A new mesh is generated using the computer code BAMG (bidimensional anisotropic mesh generator) developed
by Hecht \cite{Hec97} based on Delaunay-type triangulation. The code allows the user to supply his/her own metric
tensor defined on a background mesh.

Two typical $M_{DMP}$ and $M_{DMP+adap}$ meshes for this example are shown in Fig.~\ref{Exa5.4-1}.
The quantity $-u_{min}$ is plotted as a function
of the number of mesh elements in Fig.~\ref{Exa5.4-4}. One can see that undershoots occur
for relatively coarser meshes but not for fine meshes. This is consistent with Remark~\ref{rem4.3}
which shows that sufficiently fine $\mathbb{D}^{-1}$-uniform meshes satisfy the mesh condition (\ref{thm4.1-1}).
It should be pointed out that the maximum element size of the finest mesh for the considered range of $N$
in Fig.~\ref{Exa5.4-2} is $h \approx 0.058$, which is much larger than $h = 0.0003$ required by (\ref{rem4.3-4})
for the current example. This indicates that (\ref{rem4.3-4}) is quite conservative in estimating $h$ for a mesh
satisfying condition (\ref{thm4.1-1}). From Fig.~\ref{Exa5.4-2}, one may also notice that $M_{DMP+adap}$ meshes
lead to larger undershoots than $M_{DMP}$ meshes. Fig.~\ref{Exa5.4-5} shows that the former has a larger maximum
element size than the latter. There are no clear explanations why this should
happen. We recall that $M_{DMP}$ meshes are completely determined by $\mathbb{D}$ whereas
$M_{DMP+adap}$ meshes are determined by $\mathbb{D}$ in element shape and 
by minimization of interpolation error in element size. These two sets of meshes serve different purposes
and it seems that it can go either way.

Finally,  contours of linear finite element solutions are shown in Figs.~\ref{Exa5.4-2} and \ref{Exa5.4-3}. 
Once again, one can see that for both sets of meshes, undershoots occur for a relatively coarse mesh and
vanish for a finer mesh.
\qed
\end{exam}

\section{Conclusions}

In the previous sections
we have developed a mesh condition (\ref{thm4.1-1}) under which the linear finite element solution defined in (\ref{disc-1})
for the general anisotropic diffusion problem (\ref{bvp-pde}) and (\ref{bvp-bc}) involving convection and reaction terms
to satisfy a discrete maximum principle. 
Loosely speaking, the condition requires that the dihedral angles of all elements of the mesh
be $\mathcal{O}(\|\V{b}\|_\infty h + \|c\|_\infty h^2)-\mbox{acute}$ when they are measured in the metric specific
by the inverse of the coefficient matrix, where $\V{b}$ and $c$ are the coefficients of the convection and reaction
terms, respectively.
Moreover, we have shown that in two dimensions a weaker condition, (\ref{meshcon-4}) or (\ref{meshcon-6}) --
an $\mathcal{O}(\|\V{b}\|_\infty h + \|c\|_\infty h^2)$ perturbation of the generalized Delaunay condition
developed in \cite{Hua10}, is sufficient for the linear finite element solution to satisfy a discrete maximum principle.
Finally, it is worth pointing out that many existing mesh conditions such as those developed in
Ciarlet and Raviart \cite{CR73} (for isotropic diffusion without convection terms),
Strang and Fix \cite{SF73} (for 2D isotropic diffusion without convection terms),
Wang and Zhang \cite{WaZh11} (for isotropic diffusion with convection and reaction terms),
Li and Huang \cite{LH10} (for anisotropic diffusion without convection and reaction terms),
and Huang \cite{Hua10} (for 2D anisotropic diffusion without convection and reaction terms)
are special cases of mesh condition (\ref{thm4.1-1}) or (\ref{meshcon-4}).

\vspace{20pt}

{\bf Acknowledgment.} 
The work was supported in part by the National
Science Foundation (U.S.A.) under Grant DMS-1115118
and the Natural Science Foundation of China under Grants 10931004
and 40906048.


\begin{thebibliography}{10}

\bibitem{ABBM98a}
I.~Aavatsmark, T.~Barkve, {\O}.~B{\o}e, and T.~Mannseth.
\newblock Discretization on unstructured grids for inhomogeneous, anisotropic
  media. {I}. {D}erivation of the methods.
\newblock {\em SIAM J. Sci. Comput.}, 19:1700--1716 (electronic), 1998.

\bibitem{BKK08}
J.~Brandts, S.~Korotov, and M.~K\v{r}\'i\v{z}ek.
\newblock The discrete maximum principle for linear simplicial finite element
  approximations of a reaction-diffusion problem.
\newblock {\em Lin. Alg. Appl.}, 429:2344--2357, 2008.

\bibitem{BE04}
E.~Burman and A.~Ern.
\newblock Discrete maximum principle for {G}alerkin approximations of the
  {L}aplace operator on arbitrary meshes.
\newblock {\em C. R. Acad. Sci. Paris}, Ser.I 338:641--646, 2004.

\bibitem{Cia70}
P.~G. Ciarlet.
\newblock Discrete maximum principle for finite difference operators.
\newblock {\em Aequationes Math.}, 4:338--352, 1970.

\bibitem{Cia78}
P.~G. Ciarlet.
\newblock {\em The Finite Element Method for Elliptic Problems}.
\newblock North-Holland, Amsterdam, 1978.

\bibitem{CR73}
P.~G. Ciarlet and P.-A. Raviart.
\newblock Maximum principle and uniform convergence for the finite element
  method.
\newblock {\em Comput. Meth. Appl. Mech. Engrg.}, 2:17--31, 1973.

\bibitem{DDS04}
A.~Dr\v{a}g\v{a}nescu, T.~F. Dupont, and L.~R. Scott.
\newblock Failure of the discrete maximum principle for an elliptic finite
  element problem.
\newblock {\em Math. Comp.}, 74:1--23, 2004.

\bibitem{Eva98}
L.~C. Evans.
\newblock {\em Partial Differential Equations}.
\newblock American Mathematical Society, Providence, Rhode Island, 1998.
\newblock Graduate Studies in Mathematics, Volume 19.

\bibitem{For89}
P.~A. Forsyth.
\newblock A control-volume, finite-element method for local mesh refinement in
  thermal reservoir simulation.
\newblock {\em {SPE} Reservoir Engineering}, 5:561--566 (Paper SPE 18415),
  1990.

\bibitem{GL09}
S.~G\H{u}nter and K.~Lackner.
\newblock A mixed implicit-explicit finite difference scheme for heat transport
  in magnetised plasmas.
\newblock {\em J. Comput. Phys.}, 228:282--293, 2009.

\bibitem{GYK05}
S.~G\H{u}nter, Q.~Yu, J.~Kruger, and K.~Lackner.
\newblock Modelling of heat transport in magnetised plasmas using non-aligned
  coordinates.
\newblock {\em J. Comput. Phys.}, 209:354--370, 2005.


\bibitem{Hec97}
\newblock F. Hecht.
\newblock Bidimensional anisotropic mesh generator. Technical report, INRIA, Rocquencourt, 1997.
\newblock Source code available from: http://www.ann.jussieu.fr/hecht/ftp/bamg.



\bibitem{Hua06}
W.~Huang.
\newblock Mathematical principles of anisotropic mesh adaptation.
\newblock {\em Comm. Comput. Phys.}, 1:276-310, 2006.


\bibitem{Hua10}
W.~Huang.
\newblock Discrete maximum principle and a delaunay-type mesh condition for
  linear finite element approximations of two-dimensional anisotropic diffusion
  problems.
\newblock {\em Numer. Math. Theory Meth. Appl.}, 4:319--334, 2011.
\newblock (arXiv:1008.0562).


\bibitem{HR11}
W.~Huang and R.~D.~Russell.
\newblock {\em Adaptive Moving Mesh Methods}.
\newblock Springer, new York, 2011.


\bibitem{KK09}
J.~Kar{\'a}tson and S.~Korotov.
\newblock An algebraic discrete maximum principle in {H}ilbert space with
  applications to nonlinear cooperative elliptic systems.
\newblock {\em SIAM J. Numer. Anal.}, 47:2518--2549, 2009.

\bibitem{KKK07}
J.~Kar\'atson, S.~Korotov, and M.~K\v{r}\'i\v{z}ek.
\newblock On discrete maximum principles for nonlinear elliptic problems.
\newblock {\em Math. Comput. Sim.}, 76:99--108, 2007.

\bibitem{KSS09}
D.~Kuzmin, M.~J. Shashkov, and D.~Svyatskiy.
\newblock A constrained finite element method satisfying the discrete maximum
  principle for anisotropic diffusion problems.
\newblock {\em J. Comput. Phys.}, 228:3448--3463, 2009.

\bibitem{KL95}
M.~K\v{r}\'i\v{z}ek and Q.~Lin.
\newblock On diagonal dominance of stiffness matrices in {3D}.
\newblock {\em East-West J. Numer. Math.}, 3:59--69, 1995.

\bibitem{LePot09}
C.~Le~Potier.
\newblock A nonlinear finite volume scheme satisfying maximum and minimum
  principles for diffusion operators.
\newblock {\em Int. J. Finite Vol.}, 6:20, 2009.

\bibitem{Let92}
F.~W. Letniowski.
\newblock Three-dimensional {D}elaunay triangulations for finite element
  approximations to a second-order diffusion operator.
\newblock {\em SIAM J. Sci. Stat. Comput.}, 13:765--770, 1992.

\bibitem{LH10}
X.~P. Li and W.~Huang.
\newblock An anisotropic mesh adaptation method for the finite element solution
  of heterogeneous anisotropic diffusion problems.
\newblock {\em J. Comput. Phys.}, 229:8072--8094, 2010 (arXiv:1003.4530).

\bibitem{LSS07}
X.~P. Li, D.~Svyatskiy, and M.~Shashkov.
\newblock Mesh adaptation and discrete maximum principle for {2D} anisotropic
  diffusion problems.
\newblock Technical Report LA-UR 10-01227, Los Alamos National Laboratory, Los
  Alamos, NM, 2007.

\bibitem{LSSV07}
K.~Lipnikov, M.~Shashkov, D.~Svyatskiy, and Yu. Vassilevski.
\newblock Monotone finite volume schemes for diffusion equations on
  unstructured triangular and shape-regular polygonal meshes.
\newblock {\em J. Comput. Phys.}, 227:492--512, 2007.

\bibitem{LS08}
R.~Liska and M.~Shashkov.
\newblock Enforcing the discrete maximum principle for linear finite element
  solutions of second-order elliptic problems.
\newblock {\em Comm. Comput. Phys.}, 3:852--877, 2008.

\bibitem{MD06}
M.~J. Mlacnik and L.~J. Durlofsky.
\newblock Unstructured grid optimization for improved monotonicity of discrete
  solutions of elliptic equations with highly anisotropic coefficients.
\newblock {\em J. Comput. Phys.}, 216:337--361, 2006.

\bibitem{PM90}
P.~Perona and J.~Malik.
\newblock Scale-space and edge detection using anisotropic diffusion.
\newblock {\em {IEEE} Trans. Pattern Anal. Mach. Intel.}, 12:629--639, 1990.

\bibitem{SH07}
P.~Sharma and G.~W. Hammett.
\newblock Preserving monotonicity in anisotropic diffusion.
\newblock {\em J. Comput. Phys.}, 227:123--142, 2007.

\bibitem{Sto86}
G.~Stoyan.
\newblock On maximum principles for monotone matrices.
\newblock {\em Lin. Alg. Appl.}, 78:147--161, 1986.

\bibitem{SF73}
G.~Strang and G.~J. Fix.
\newblock {\em An Analysis of the Finite Element Method}.
\newblock Prentice Hall, Englewood Cliffs, NJ, 1973.

\bibitem{WaZh11}
J.~Wang and R.~Zhang.
\newblock Maximum principle for {P1}-conforming finite element approximations
  of quasi-linear second order elliptic equations.
\newblock {\em SIAM J. Numer. Anal.}, 50:626--642, 2012.
\newblock (arXiv:1105.1466).

\bibitem{Wei98}
J.~Weickert.
\newblock {\em Anisotropic Diffusion in Image Processing}.
\newblock Teubner-Verlag, Stuttgart, Germany, 1998.

\bibitem{XZ99}
J.~Xu and L.~Zikatanov.
\newblock A monotone finite element scheme for convection-diffusion equations.
\newblock {\em Math. Comput.}, 69:1429--1446, 1999.


\end{thebibliography}

\end{document}